\documentclass[a4paper, 12pt]{article}
\usepackage{enumerate, theorem}
\usepackage{amsmath, amsfonts, amssymb}
\usepackage[height=22.5cm, width=15cm]{geometry}
\usepackage[all]{xy}
\usepackage{mathrsfs, yfonts}

\newcommand{\A}{\tilde{\mathcal{A}}}

\newcommand{\parag}[1]{\paragraph{\sc{#1.}}}

\newtheorem{thm}{Th\'eor\`eme}[subsection]
\newtheorem{defn}[thm]{D\'efinition}
\newtheorem{cor}[thm]{Corollaire}
\newtheorem{prop}[thm]{Proposition}
\newtheorem{lemma}[thm]{Lemme}

\begin{document}

\title{P\'eriodes \'evanescentes et  (a,b)-modules monog\`enes.}

\date{05/01/09}

\author{Daniel Barlet\footnote{Barlet Daniel, Institut Elie Cartan UMR 7502  \newline
Nancy-Universit\'e, CNRS, INRIA  et  Institut Universitaire de France, \newline
BP 239 - F - 54506 Vandoeuvre-l\`es-Nancy Cedex.France.\newline
e-mail : barlet@iecn.u-nancy.fr}.}

 \maketitle
 
 \markright{P\'eriodes \'evanescentes ...}

 \section*{Abstract.}
 In order to describe the asymptotic behaviour of a  vanishing period in a one parameter family we introduce and use  a very simple algebraic structure : regular geometric (a,b)-modules generated (as left \ $\A-$modules) by one element. The idea is to use not the full Brieskorn module associated to the Gauss-Manin connection but a minimal (regular) differential equation satisfied by the period integral we are interested in. We show that the Bernstein polynomial associated is quite simple to compute for such (a,b)-modules and give a precise description of the exponents which appears in the asymptotic expansion which avoids integral shifts. We show a couple of explicit computations in some classical (but not so easy) examples. 
 
 \bigskip

AMS Classification (2000) : 32-S-25, 32-S-40, 32-S-50.

\bigskip
 
 Key words : vanishing cycles, vanishing periods, (a,b)-modules, Brieskorn modules, Gauss-Manin systems.
 
 \newpage

  \tableofcontents

 \section{Introduction}
 
Cet article a pour but d'\'etudier le d\'eveloppement asymptotique d'une int\'egrale de la forme
$$ s \to \int_{\gamma_s} \frac{d\omega}{df} $$
issue d'une des deux situations g\'eom\'etriques suivantes :
\begin{enumerate}
\item Soit \ $f : X \to D$ \ un repr\'esentant de Milnor d'un germe de fonction holomorphe \`a l'origine de \ $\mathbb{C}^{n+1}$, soit \ $\omega$ \ une \ $p-$forme holomorphe sur \ $X$ \ v\'erifiant \ $df\wedge d\omega \equiv 0 $, et soit \ $(\gamma_s)_{s \in D}$ \ une famille horizontale de \ $p-$cycles compacts dans les fibres de \ $f$.
\item Soit \ $f : X \to D$ \ une fonction holomorphe propre d'une vari\'et\'e complexe \ $X$ \ de dimension \ $n+1$ \ dont les singularit\'es  \'eventuelles sont toutes contenues dans la fibre \ $Y : = f^{-1}(0)$, soit \ $\omega$ \ une \ $p-$forme holomorphe sur \ $X$ \ v\'erifiant \ $df\wedge d\omega \equiv 0 $, et soit \ $(\gamma_s)_{s \in D}$ \ une famille horizontale de \ $p-$cycles compacts dans les fibres de \ $f$.
\end{enumerate}
La connexion de Gauss-Manin de \ $f$ \ en degr\'e \ $p$ \ est un syst\`eme diff\'erentiel sur \ $D$ \ avec une singularit\'e r\'eguli\`ere \`a l'origine qui "contr\^ole" ce type de fonction quand la forme diff\'erentielle \ $\omega$ \ varie. Mais il est clair que pour une \ $p-$forme donn\'ee, le "contr\^ole" en question est assez impr\'ecis. \\
En effet l'int\'egrale consid\'er\'ee est solution d'une \'equation diff\'erentielle "minimale" (en un sens \`a pr\'eciser) qui est "contenue" dans le syst\`eme de Gauss-Manin, mais peut-\^etre beaucoup plus "petite", et donne donc des renseignements beaucoup plus pr\'ecis sur la fonction consid\'er\'ee et les termes int\'eressants apparaissant dans son d\'eveloppement asymptotique (convergeant en fait).

\smallskip

Nous nous int\'eresserons plus pr\'ecis\'ement aux cas o\`u l'on sait d\'ej\`a que le syst\`eme de Gauss-Manin consid\'er\'e est incarn\'e par un (a,b)-module r\'egulier\footnote{Le concept de (a,b)-module g\'en\'eralise le compl\'et\'e formel d'un module de Brieskorn d'une singularit\'e isol\'ee ; voir par exemple [S. 89], [B. 93] ou [B. 05].}. Le cas classique est \'evidemment celui d'un germe \ $f$ \ \`a singularit\'e isol\'ee \`a l'origine (que l'on peut aussi voir comme le cas global \`a singularit\'e isol\'ee). Pour le cas d'un germe quelconque, quitte \`a quotienter par la torsion, on est toujours dans cette situation d'apr\`es [B.-S. 04]. Mais la pr\'esence de  torsion non de type fini peut-\^etre une difficult\'e s\'erieuse pour mener \`a bien des calculs concrets. Plus r\'ecemment on a montr\'e dans [B.II] compl\'et\'e par [B.07] que pour un germe \ $f$ \ \`a singularit\'e de dimension 1 la torsion est de type finie. On a aussi montr\'e dans [B.08]  que c'est toujours la situation pour une fonction holomorphe propre (situation du 2.  ci-dessus).

\smallskip

Dans le cadre des (a,b)-modules (ou modules de Brieskorn g\'en\'eralis\'es) il s'agit de consid\'erer un (a,b)-module r\'egulier g\'eom\'etrique \ $E$\footnote{La r\'egularit\'e du (a,b)-module correspond ici \`a la r\'egularit\'e de la connexion de Gauss-Manin, l'aspect g\'eom\'etrique encode \`a la fois le th\'eor\`eme de Monodromie et le th\'eor\`eme de positivit\'e de B. Malgrange ; voir [M. 75].} et de consid\'erer le sous-(a,b)-module minimal contenant un \'el\'ement \ $x \in E$ \ donn\'e. Ceci conduit naturellement \`a la notion de (a,b)-module {\bf monog\`ene}  qui est introduite et d\'ecrite dans cet article. Un invariant important est l'\'el\'ement de Bernstein d'un tel (a,b)-module monog\`ene qui est dans ce cas un avatar tr\`es simple du polyn\^ome de Bernstein d'un (a,b)-module r\'egulier (voir [B.93]). Le calcul de cet invariant, qui est beaucoup plus simple que le calcul du (a,b)-module engendr\'e par \ $x$ \ lui-m\^eme, c'est \`a dire la d\'etermination pr\'ecise de l'\'equation diff\'erentielle "minimale" v\'erifi\'ee par la fonction consid\'er\'ee, donne d\'ej\`a des renseignements pr\'ecis sur les exposants apparaissant dans le d\'eveloppement asymptotique \`a l'origine de l'int\'egrale consid\'er\'ee.

\smallskip

De mani\`ere \`a bien montrer l'aspect "concret" de notre approche, nous avons d\'etaill\'e deux exemples dans le cas classique d'une singularit\'e isol\'ee, qui font bien appara\^itre le gain en pr\'ecision par rapport \`a ce que donne le calcul du syst\`eme de Gauss-Manin complet dans ces cas, ou \`a fortiori par rapport \`a un calcul topologique de la monodromie.

\bigskip

D\'ecrivons de fa{\c c}on un peu plus pr\'ecise les principaux r\'esultats obtenus.\\

Nous commen{\c c}ons par montrer un r\'esultat  de factorisation (proposition  \ref{homogene}) des \'el\'ements homog\`enes en (a,b) dans l'alg\`ebre \ $\A$ \ qui met en \'evidence une repr\'esentation "tordue" du groupe sym\'etrique. Ceci joue un r\^ole important dans la compr\'ehension du calcul du polyn\^ome de Bernstein et les th\'eor\`emes de structure des (a,b)-modules monog\`enes r\'eguliers.

\smallskip

Le premier th\'eor\`eme de structure \ref{Structure} met en \'evidence le lien entre la forme initiale du g\'en\'erateur de l'id\'eal annulateur d'un g\'en\'erateur d'un (a,b)-module monog\`ene r\'egulier et son polyn\^ome de Bernstein. L'aspect multiplicatif de l'\'el\'ement de Bernstein associ\'e en est une cons\'equence importante. Ce th\'eor\`eme met \'egalement en \'evidence des invariants  plus fins  (voir \ref{Invariants}) qui seront exploit\'es dans un travail ult\'erieur.

\smallskip

Le second th\'eor\`eme de structure \ref{Factorisation 1} est multiplicatif et il explicite le lien entre les suites de Jordan-H{\"o}lder et l'id\'eal annulateur d'un g\'en\'erateur d'un (a,b)-module monog\`ene r\'egulier.\\
Nous donnons ensuite un th\'eor\`eme (voir  \ref{Utile}) permettant de d\'eduire la valeur du polyn\^ome de Bernstein d'un (a,b)-module monog\`ene r\'egulier d\`es que l'on connait un \'el\'ement convenable annulant son g\'en\'erateur sans pour autant demander de conna\^itre la forme canonique donn\'ee dans le premier th\'eor\`eme de structure. Ceci est l'outil crucial qui permet de mener a bien les calculs d'exemples qui concluent cet article.\\
Enfin nous montrons (voir le th\'eor\`eme \ref{Realisation}) que tout (a,b)-module monog\`ene g\'eom\'etrique est engendr\'e par une fonction multiforme de d\'etermination finie "de type standard" et r\'eciproquement. Ceci permet de montrer que la connaissance du polyn\^ome de Bernstein r\'epond \`a la question "basique" de la d\'etermination {\em pr\'ecise} (sans d\'ecalage entier et avec l'exposant du logarithme) des termes "fondamentaux" du d\'evelopement asymptotique standard consid\'er\'e.

  \section{Polyn\^omes homog\`enes en (a,b), unitaires en  a.}
  
  On notera par  \ $ \A $ \ la \ $\mathbb{C}[[b]]-$alg\`ebre
  $$\A : = \{ \sum_{\nu = 0}^{\infty} \ P_{\nu}(a).b^{\nu}, P_{\nu} \in \mathbb{C}[z]\} $$
 o\`u l'on a \ $ a.b - b.a = b^2$, et o\`u la multiplication par \ $a$ \ \`a gauche ( ou \`a droite) est continue pour la topologie \ $b-$adique.
 
 \begin{defn} Nous dirons qu'un \'el\'ement de \ $\A$ \ est {\bf homog\`ene de degr\'e \ $k$ \ en \ $(a,b)$} \ si on peut l'\'ecrire comme une combinaison lin\'eaire
 $$ \sum_{j=0}^k \ \lambda_j.a^j.b^{k-j}  $$
 o\`u les \ $\lambda_j$ \ sont des nombres complexes. On dira qu'un tel \'el\'ement est {\bf unitaire en $a$} si on a \ $\lambda_k = 1$.
 \end{defn}
 
 Le lecteur v\'erifiera facilement que tout \'el\'ement de \ $\A$ \ de la forme \ $$a^{m_1}b^{n_1}a^{m_2}b^{m_2}\cdots a^{m_h}b^{n_h} $$
 avec \ $ \ \sum_{p=1}^h \ m_p + n_p = k$, est homog\`ene de degr\'e \ $k$ \ en \ $(a,b)$. Ceci montre que le produit d'un \'el\'ement homog\`ene de degr\'e \ $k$ \ par un \'el\'ement homog\`ene de degr\'e \ $l$ \ est homog\`ene de degr\'e \ $k+l$.\\
  
  On remarquera que les mon\^omes \ $b^j.a^{k-j}, j \in [0,k]$ \ forment \'egalement une base de l'espace vectoriel des \'el\'ements homog\`enes de degr\'e \ $k$ \ de l'alg\`ebre \ $\A$.
 
 \bigskip
 
 Notre premier objectif est de montrer la proposition suivante :
 
 \begin{prop}\label{homogene}
 Tout \'el\'ement homog\`ene \ $P \in \A$ \ de degr\'e \ $k$ \ en \ $(a,b)$ \  et unitaire en \ $a$ \ peut s'\'ecrire sous la forme \ $ P =  (a - \lambda_1.b)(a - \lambda_2.b)\cdots (a - \lambda_k.b) $, o\`u les \ $\lambda_j$ \ sont des nombres complexes. De plus, pour chaque \ $\sigma \in \mathfrak{S}_k$ \ on aura encore
 $$ P = (a-(\lambda_{\sigma(1)}+\sigma(1)-1).b)\cdots (a-(\lambda_{\sigma(k)}+\sigma(k)-k).b) $$
 et on obtient ainsi, quand \ $\sigma$ \ d\'ecrit le groupe \ $\mathfrak{S}_k$ \ des permutations de \ $\{1,2,\cdots,k\}$ \ toutes les fa{\c c}ons de d\'ecomposer \ $P$ \ en produit d'\'el\'ements homog\`enes de degr\'e $1$ \ en \ $(a,b)$.
 \end{prop}
 
 \parag{Preuve} Commen{\c c}ons par remarquer que l'application
 $$ \varphi : \mathfrak{S}_k \times \mathbb{C}^k \to \mathbb{C}^k $$
 d\'efinie par \ $\varphi(\sigma,\lambda_1, \cdots, \lambda_k) = (\lambda_{\sigma(1)}+\sigma(1)-1, \cdots, \lambda_{\sigma(k)}+\sigma(k)-k)$ \ est une action du groupe \ $\mathfrak{S}_k$ \ sur \ $\mathbb{C}^k$. En effet on a
 \begin{align*}
 &  \varphi(\tau,\varphi(\sigma,(\lambda_j))) = \varphi(\tau,(\lambda_{\sigma(j)}+\sigma(j)-j)) \\
 & \qquad   = (\lambda_{\tau(\sigma(j))} + \sigma(j) - j + \tau(\sigma(j)) - \sigma(j)) = ( \lambda_{\tau(\sigma(j))} + \tau(\sigma(j)) - j) \\
 & \qquad   = \varphi(\tau\circ\sigma, (\lambda_j)).
 \end{align*}
 Nous appellerons cette action {\bf"l'action tordue"} de \ $\mathfrak{S}_k$ \ sur \ $\mathbb{C}^k$.\\
 Pour montrer l'\'egalit\'e 
 $$(a - \lambda_1.b)(a - \lambda_2.b)\cdots (a - \lambda_k.b)  = (a - (\lambda_{\sigma(1)}+\sigma(1)-1).b)\cdots (a - (\lambda_{\sigma(k)}+\sigma(k)-k).b) $$
 il suffit donc de traiter le cas de la transposition de \ $\mathfrak{S}_2$, c'est-\`a-dire de v\'erifier dans \ $\A$ \ l'\'egalit\'e \ $ (a - \lambda_1.b).(a - \lambda_2.b) = (a - (\lambda_2 +1).b).(a - (\lambda_1- 1).b)$, ce qui est imm\'ediat.\\
 Passons \`a l'existence de la d\'ecomposition de \ $P$ \ en produit d'\'el\'ements homog\`enes de degr\'e \ $1$ \ en \ $(a,b)$. Pour cela montrons que \ $ b^{-k}.P$ \ est un polyn\^ome unitaire de degr\'e \ $k$ \ en \ $b^{-1}.a$. Comme il suffit de prouver que pour chaque \ $j \in [0,k]$ \ l'\'el\'ement \ $b^{-k}.b^j.a^{k-j} \in \A[b^{-1}]$ \ est un polyn\^ome en \ $b^{-1}.a$ \ de degr\'e au plus \ $k$, \'ecrivons, puisque l'on a \ $ a.b^m = b^m.a + m b^{m+1} \quad \forall m \in \mathbb{Z}$,
 \begin{align*}
 &  b^{-k}.b^j.a^{k-j} =b^{-1}.(b^{j-k+1}.a).a^{k-j-1} \\
 & \qquad \qquad = b^{-1}.(a.b^{j-k+1} - (j-k+1)b^{j-k+2}).a^{k-j-1}\\
 & \qquad \qquad = b^{-1}.a.(b^{j-k+1}.a^{k-j-1}) - (j-k+1)b^{j-k+1}.a^{k-j-1} \\
 & \qquad \qquad = (b^{-1}.a - (j-k+1)).b^{-k}.b^{j+1}.a^{k-j-1} 
 \end{align*}
 ce qui prouve notre assertion par r\'ecurrence descendante sur \ $j$, le cas \ $j = k$ \ \'etant clair.\\
 La factorisation du polyn\^ome unitaire de degr\'e \ $k$ \ ainsi obtenu donne alors la factorisation d\'esir\'ee pour \ $P$,  de la fa{\c c}on suivante : \\
 La relation
 $$ b^k.\prod_{j=1}^k \ (b^{-1}.a - \lambda_j) = (a - (\lambda_1+k -1).b)\cdots (a - (\lambda_k + k - k).b) $$
 s'obtient facilement en remarquant que l'on a
 $$ b^k.(b^{-1}.a - \lambda) = b^{k-1}.(a - \lambda.b) = (a - (\lambda +k-1).b).b^{k-1} $$
 ce qui permet de conclure. $\hfill \blacksquare$
 
 \bigskip
 
 Ceci montre la correspondance entre les orbites de l'action usuelle de \ $\mathfrak{S}_k$ \ sur \ $\mathbb{C}^k$ \ et les orbites de l'action tordue via la correspondance qui \`a un polyn\^ome unitaire \ $Q$ \ de degr\'e \ $k$ \ associe l'\'el\'ement homog\`ene de degr\'e \ $k$ \ en \ $(a,b)$, unitaire en \ $a$ \ d\'efini par \ $ P : = b^k.Q(b^{-1}.a)$. On a en fait un isomorphisme 
 $$ \varphi : \mathbb{C}^k \to  \mathbb{C}^k $$
 donn\'e par \ $\varphi(\lambda_1, \cdots, \lambda_k) = (\lambda_1+k-1,\cdots, \lambda_k +k - k)$ \ qui est \'equivariant quand on fait agir \ $\mathfrak{S}_k$ \ au d\'epart par l'action usuelle et \`a l'arriv\'ee par l'action tordue.

 \bigskip

 La proposition qui suit donne  une autre fa{\c c}on de voir cette action tordue et l'\'equivariance d\'ecrite ci-dessus. 
 
 \begin{prop}\label{p}
 Soient \ $\mu_1, \cdots, \mu_k$ \ des nombres complexes deux \`a deux distincts. Alors on a l'\'egalit\'e suivante entre id\'eaux \`a gauche de l'alg\`ebre \ $\A$ :
 \begin{equation*} \cap_{j=1}^k \ \A.(a - \mu_j.b) = \A.(a - (\mu_1+k-1).b)\cdots (a - (\mu_k+k-k).b) . \tag{*}
 \end{equation*}
 \end{prop}
 
 La d\'emonstration utilisera le lemme suivant :

 \begin{lemma}\label{2}
 Soient \ $\lambda \not= \mu$ \ deux nombres complexes distincts et soit \ $x \in \A$. Alors \ $x.(a - \lambda.b) \in \A.(a - \mu.b) $ \ si et seulement si \ $x \in \A.(a - (\mu+1).b)$.
 \end{lemma}
 
 \parag{Preuve} En utilisant l'automorphisme unitaire de la \ $\mathbb{C}-$alg\`ebre  \ $\A$ \ qui envoie \ $a$ \ sur \ $a - \mu.b$ \ et \ $b$ \ sur \ $b$, on se ram\`ene au cas o\`u  \ $\lambda \not= \mu = 0$.\\
  L'\'egalit\'e \ $ (a - b).(a -\lambda.b) = (a - (\lambda-1).b).a $ \ montre que la condition \ $x \in \A.(a-b)$ \ implique bien \ $x.(a - \lambda.b) \in \A.a $. \\
 R\'eciproquement, supposons que \ $x.(a-\lambda.b) = t.a $. \'Ecrivons \ $x = z.(a-b) + r $ \ avec \ $r \in \mathbb{C}[[b]]$. On obtient alors
 $$ x.(a - \lambda.b) = z.(a-b).(a - \lambda.b) + r.(a - \lambda.b)= t.a $$
 ce qui donne \ $\lambda.r.b \in \A.a $. Comme \ $\lambda \not= 0$ \ on en conclut que \ $r = 0$ \ ce qui prouve l'assertion \ $x \in \A.(a - b)$. $\hfill \blacksquare$
 
 \parag{Preuve de la proposition \ref{p}} L'inclusion \ $\supseteq$ \ de l'\'enonc\'e r\'esulte imm\'ediatement de l'invariance du produit consid\'er\'e dans le second membre de \ $(^*)$ \  par l'action tordue de \ $ \mathfrak{S}_k$. Montrons l'inclusion oppos\'ee par r\'ecurrence sur \ $k$. Le cas \ $k = 1$ \ \'etant clair, supposons l'assertion montr\'ee pour \ $k-1 \geq 1$.\\
  Soit donc \ $x \in  \cap_{j=1}^k \big( \A.(a - \mu_j.b)\big)$. \'Ecrivons \ $ x = y.(a-\mu_k.b)$. On a alors pour chaque \ $j \in [1, k-1] \quad y.(a-\mu_k.b) \in \A.(a - \mu_j.b)$. \\
De plus l'hypoth\`ese que les \ $\mu_j$ \ sont deux \`a deux distincts pemet d'appliquer le lemme \ref{2} pour chaque \ $j \in [1,k-1]$, ce qui donne que \ $y \in \cap_{j=1}^{k-1} \A.(a - \mu_j.b) $. Alors l'hypoth\`ese de r\'ecurrence permet de conclure. $\hfill \blacksquare$ 

\bigskip

Terminons ce paragraphe en montrant que tout (a,b)-module r\'egulier contient un  sous-module monog\`ene qui est de codimension finie.

\begin{lemma}
Pour tout (a,b)-module r\'egulier \ $E$ \ il existe \ $x \in E$ \ tel que le sous-(a,b)-module monog\`ene \ $\A.x$ \ de \ $E$ \ soit de codimension finie dans \ $E$.
\end{lemma}

\parag{Preuve} Montrons ce r\'esultat par r\'ecurrence sur \ $k$, le rang de \ $E$. Le cas du rang 1 \ \'etant \'evident, supposons l'assertion montr\'ee en rang \ $k-1$, avec \ $k \geq 2$ \ et montrons-l\`a au rang \ $k$. Consid\'erons une suite exacte de (a,b)-modules\footnote{une telle suite existe d'apr\`es [B.93] prop.2.2.}
$$ 0 \to E_{\lambda} \to E \overset{\pi}{\to} F \to  0  $$
o\`u \ $E$ \ est de rang \ $k$, et soit \ $y \in F$ \ tel que \ $\A.y$ \ soit de codimension finie dans \ $F$. Soit \ $x_0 \in E$ \ v\'erifiant \ $ \pi(x_0) = y$. Nous allons consid\'erer deux cas.
\parag{Premier cas : $ \A.x_0 \cap E_{\lambda} \not= \{0\}$} Il existe alors un entier \ $\nu \geq 0 $ \ tel que \ $\A.x_0 \cap E_{\lambda} = b^{\nu}.E_{\lambda}  $. On a alors la suite exacte de (a,b)-modules
$$ 0 \to b^{\nu}.E_{\lambda} \to \A.x_0 \to \A.y \to 0 $$
qui montre que le rang de \ $\A.x_0$ \ est \'egale \`a \ $k$. Ce qui permet de conclure.

\parag{Second cas : $ \A.x_0 \cap E_{\lambda} = \{0\}$} Soit \ $P \in \A$ \ le polyn\^ome unitaire en \ $a$ \ qui engendre l'annulateur \`a gauche dans \ $\A$ \ de l'\'el\'ement \ $x_0 \in E$. Montrons qu'il existe un entier \ $\nu \geq 0$ \ tel que \ $P.b^{\nu}.e_{\lambda} \not= 0 $ \ o\`u \ $e_{\lambda}$ \ d\'enote un g\'en\'erateur standard de \ $E_{\lambda}$ \ (donc v\'erifiant \ $a.e_{\lambda} = \lambda.b.e_{\lambda}$). Chaque entier \ $\nu$ \ tel que \ $P.b^{\nu}.e_{\lambda} = 0 $ \ d\'efinit un morphisme (a,b)-lin\'eaire
$$ \varphi_{\nu} : \A\big/ \A.P \to  E_{\lambda}$$
en posant \ $\varphi_{\nu}(1) : = b^{\nu}.e_{\lambda}$. De plus, pour des entiers \ $\nu_1, \cdots, \nu_N$ \ deux \`a deux distincts, on obtient ainsi des \'el\'ements \ $\mathbb{C}-$lin\'eairement ind\'ependants dans \ $Hom_{\A}(\A\big/ \A.P, E_{\lambda})$ \ puisque l'image de \ $\varphi_{\nu}$ \ est \'egale \`a \ $b^{\nu}.E_{\lambda}$. Comme l'espace vectoriel \ $Hom_{\A}(\A\big/ \A.P, E_{\lambda})$ \ est de dimension finie d'apr\`es [B.95] (corollaire du th. 1 ter), on en conclut qu'il n'existe qu'un nombre fini d'entiers \ $\nu$ \ v\'erifiant \ $P.b^{\nu}.e_{\lambda} = 0 $, ce qui prouve notre assertion.\\
Fixons alors \ $\nu$ \ tel que \ $P.b^{\nu}.e_{\lambda} \not= 0 $ \ et posons \ $x_1 : = x_0 + b^{\nu}.e_{\lambda}$. On a alors \ $\pi(x_1) = y $ \ et \ $\A.x_1 \cap E_{\lambda} \not= \{0\}$ \ puisque \ $P.x_1 = P.b^{\nu}.e_{\lambda} \in E_{\lambda} \setminus \{0\}$. Nous sommes donc ramen\'e au premier cas, quitte \`a choisir \ $x_1$ \ \`a la place de \ $x_0$. $\hfill \blacksquare$

\section{Les th\'eor\`emes de structure.}

\subsection{Pr\'eambule.}

\begin{defn}
On dira qu'un (a,b)-module \ $E$ \ est {\bf monog\`ene}  s'il est isomorphe comme \ $\A-$module \`a gauche \`a un quotient \ $\A\big/I$ \ o\`u \ $I$ \ est un id\'eal \`a gauche de \ $\A$. $\hfill \square$
\end{defn}

\bigskip

Rappelons qu'un (a,b)-module est, par d\'efinition un \ $\A-$module \`a gauche qui est libre et de type fini sur la sous-alg\`ebre \ $\mathbb{C}[[b]]$. Donc tout id\'eal \`a gauche de \ $\A$ \ ne donnera pas un (a,b)-module, mais quand c'est le cas il est n\'ecessairement monog\`ene. Le th\'eor\`eme \ref{Structure} caract\'erisera les id\'eaux qui conviennent.

\parag{Exemple} Soit \ $E$ \ un (a,b)-module et soit \ $e \in E$. Alors \ $\A.e \subset E$ \ est un (a,b)-module monog\`ene, puisque \ $\mathbb{C}[[b]]$ \ est noeth\'erien et \ $E$ \ sans \ $b-$torsion. Si \ $I$ \ est l'annulateur de \ $e$ \ dans \ $E$, on a un isomorphisme "\'evident" de \ $\A-$modules \`a gauche\ $ \varphi : \A\big/I \to \A.e $ \ d\'efini en posant \ $\varphi(1) = e$.

\bigskip

 \begin{lemma}\label{petit}
   Soit \ $E$ \ un (a,b)-module  de rang \ $k$. On suppose que \ $E$ \ est local, c'est \`a dire qu'il existe \ $N \in \mathbb{N}$ \ tel que \ $a^N.E \subset b.E$. Alors les propri\'et\'es suivantes pour un \'el\'ement \ $x \in E$ \ sont \'equivalentes
   \begin{enumerate}[1)]
   \item \ $x$ \ est un g\'en\'erateur de \ $E$ \ comme \ $\A-$module.
   \item \ $\{x, a.x, \cdots, a^{k-1}.x\}$ \ induit une \ $\mathbb{C}-$base de \ $E/b.E$.
   \item  \ $\{x, a.x, \cdots, a^{k-1}.x\}$ \ est \ $\mathbb{C}[[b]]-$base de \ $E$.
   \item \ $x$ \ engendre le \ $\mathbb{C}-$espace vectoriel \ $E/a.E + b.E $.
   \end{enumerate}
   \end{lemma} 
   
   \parag{Preuve} Montrons que \ $1) \Leftrightarrow 2) $. Si \ $x$ \ est un g\'en\'erateur de \ $E$, il engendre \ $E/b.E$ \ comme \ $\mathbb{C}[a]-$module. Comme \ $E$ \ est local, l'action de \ $a$ \ sur \ $E/b.E$ \ est nilpotente, donc \ $2)$ \ est v\'erifi\'ee. R\'eciproquement supposons \ $2)$ \ v\'erifi\'ee. Pour \ $y \in E$ \ construisons par r\'ecurrence une suite de polyn\^omes \ $P_n$ \ de degr\'es \ $\leq k-1$ \ tels que l'on ait \ $y - \sum_{j=0}^n \ b^j.P_j(a).x \in b^{n+1}.E$. Pour construire \ $P_{n+1}$ \ posons \ $y - \sum_{j=0}^n \ b^j.P_j(a).x =  b^{n+1}.z_n $ \ et \'ecrivons, gr\^ace \`a 2)
   $$ z_n = \sum_{h=0}^{k-1} \lambda_h.a^h.x + b.z_{n+1}.$$
   Posons \ $P_{n+1}(a) : = \sum_{h=0}^{k-1} \lambda_h.a^h$. On aura alors
   $$y - \sum_{j=0}^{n+1} \ b^j.P_j(a).x = b^{n+2}.z_{n+1} .$$
   En posant alors \ $ u : = \sum_{j=0}^{\infty} \ b^j.P_j(a)\in \A$, on obtient \ $ y = u.x$ \ puisque \ $E$ \ est complet pour la topologie \ $b-$adique.\\
   Les assertions \ $2) \Rightarrow 3)$ \ et \ $3) \Rightarrow 4)$ \ sont faciles et laiss\'ees au lecteur.\\
   Montrons \ $4) \Rightarrow 2)$. Soit \ $k$ \ minimal tel que \ $a^k.E \subset b.E$.  Si \ $y \in E$ \ construisons par r\'ecurrence les nombres complexes \ $\lambda_1, \cdots, \lambda_{k-1}$ \ tels que  $$y - \sum_{j=0}^h \lambda_j.a^j.x \in a^{h+1}.E + b.E .$$
    Si \ $\lambda_1, \cdots, \lambda_h$ \ sont construits, posons
   $$ y - \sum_{j=0}^h \lambda_j.a^j.x =a^{h+1}z + b.\xi $$
   et \'ecrivons \ $z = \lambda_{h+1}.x + a.t + b.v $ \ gr\^ace \`a notre hypoth\`ese. Alors on aura
   $$  y - \sum_{j=0}^{h+1} \lambda_j.a^j.x \in a^{h+2}.E + b.E $$
   ce qui fait avancer la r\'ecurrence.\\
   On aura donc \ $ y - \sum_{j=0}^{k-1} \lambda_j.a^j.x \in b.E$, ce qui montre que \ $x, a.x, \cdots, a^{k-1}.x$ \ est un syst\`eme g\'en\'erateur de \ $E/b.E$. On a donc prouv\'e \ $2)$.$\hfill \blacksquare$
   
   \parag{Remarque} Le fait qu'il existe \ $x \in E$ \ v\'erifiant la condition 1), c'est \`a dire le fait d'\^etre monog\`ene, \'equivaut donc, pour un (a,b)-module local, au fait que l'espace vectoriel \ $E/a.E + b.E$ \ soit de dimension 1, c'est \`a dire \`a l'existence d'un \ $x \in E$ \ v\'erifiant la condition 4). $\hfill \square$
   
   \bigskip
   
   \begin{cor} 
   Soit \ $E$ \ un (a,b)-module monog\`ene local de rang \ $k$ \  et   soit \ $F$ \ un sous-(a,b)-module normal de rang \ $k-1$ \ de \ $E$;   soit \ $x$ \  un g\'en\'erateur de \ $E$, \ $\lambda $ \ un nombre complexe et soit  \ $\xi \in a^2.E + a.b.E + b^2.E $. Si \ $y : =(a - \lambda.b).x + \xi  \in F$ \ alors \ $y$ \ est un g\'en\'erateur de \ $F$. 
   \end{cor} 
   
     \parag{Preuve}
   Comme \ $F$ \ est normal dans \ $E$, l'application \'evidente \ $F/b.F \to E/b.E$ \ est injective. Donc le fait que \ $a$ \ soit nilpotent sur \ $E/b.E$ \ impose la m\^eme condition sur \ $F/b.F$, ce qui montre que \ $F$ \ est local. \\
   Pour prouver que \ $F$ \ est monog\`ene, 
   remarquons que, comme \ $F$ \ est normal dans \ $E$, $F/b.F$ \ est un hyperplan de \ $E/b.E$. Comme on a \ $ y - a.x \in a^2.E + b.E$, on obtient pour chaque \ $h \in [0,k-2]$ :
   $$ a^h.y = a^{h+1}.x + \sum_{j=2}^{k-h-1} \lambda_{j,h}.a^{j+h}.x  \quad {\rm modulo} \ b.E $$
   ce qui montre que \ $y, a.y, \cdots, a^{k-2}.y$ \ est un syst\`eme libre de \ $F/b.F$ \ et donc une base. On conclut gr\^ace \`a l'implication \ $2) \Rightarrow 1)$ \ du lemme pr\'ec\'edent appliqu\'ee \`a \ $F$. $\hfill \blacksquare$
   
\bigskip

\begin{prop}[Isomorphisme]\label{Iso}
Soit \ $E$ \ un (a,b)-module monog\`ene r\'egulier et soient \ $x$ \ et \ $y$ \ deux g\'en\'erateurs de \ $E$ \ d'annulateurs respectifs \ $I$ \ et \ $J$ \ dans \ $\A$. Alors il existe \ $u,v \in \A$ \ v\'erifiant les propri\'et\'es suivantes :
\begin{align*}
&   1 - vu \in J \quad {\rm et} \quad  \ I = \{ \beta \in \A \ \big/ \  \beta.u \in J  \} . \tag{@}  \\
&  1 - uv \in I  \quad  {\rm et} \quad  \  J = \{ \alpha \in \A \ \big/ \  \alpha.v \in I \}.  \tag{@'}
\end{align*}
R\'eciproquement, si \ $\A\big/ I$ \ et \ $\A\big/ J$ \ sont des (a,b)-modules monog\`enes r\'eguliers isomorphes s'il existe \ $u,v \in \A$ \ v\'erifiant la condition \ $(@)$.
\end{prop}

\parag{Preuve} Comme \ $x$ \ et \ $y$ \ sont des g\'en\'erateurs de \ $E$, il existe \ $u,v \in \A$ \ tels que l'on ait \ $ x = u.y$ \ et \ $y = v.x$ \ dans \ $E$. On en d\'eduit que \ $1 - uv \in I$ \ et que \ $1 - vu \in J$. De plus on a pour \ $\beta \in \A$ \ l'\'equivalence entre \ $\beta.x = 0$ \ et \ $\beta.u.y = 0$ \ c'est \`a dire la seconde \'egalit\'e de \ $(@)$. Celle de \ $(@')$ \ est analogue.\\
R\'eciproquement, si l'on suppose que \ $(@)$ \ est v\'erifi\'ee, d\'efinissons l'application \ $\A-$lin\'eaire \`a gauche
$$ \varphi : \A \to \A\big/ J $$
en posant \ $\varphi(1) : = u $. Alors le noyau de \ $\varphi$ \ est l'ensemble des \ $\beta \in \A$ \ tels que \ $\beta.u \in J$. C'est donc exactement \ $I$. Donc \ $\varphi$ \ induit une injection \ $\A-$lin\'eaire  \ $\psi : \A\big/ I \to \A\big/ J$. Mais comme \ $\varphi(v) = vu = 1 \ modulo \  J$, l'application \ $\varphi$ \ est surjective et \ $\psi$ \ est un isomorphisme de \ $\A-$modules. $\hfill \blacksquare$

\parag{Remarque} On montrera plus loin (voir l'exemple qui suit  le corollaire \ref{}) qu'il est possible que \ $I$ \ et \ $J$ \ v\'erifient \ $(@)$ \ sans que les \'el\'ements \ $u$ \ et \ $v$ \  soient des inversibles de \ $\A$. $\hfill \square$

\bigskip

\subsection{Le premier  th\'eor\`eme de structure.}

\bigskip 

\begin{thm}\label{Structure}
Soit \ $I$ \ un id\'eal \`a gauche de \ $\A$. Alors le quotient \ $\A\big/I$ \ est un (a,b)-module r\'egulier de rang \ $k$ \ si et seulement si l'id\'eal \ $I$ \ est principal et admet un g\'en\'erateur de la forme
$$  x =  P + b^2.R $$
o\`u \ $P \in \A$ \ est un \'el\'ement homog\`ene de degr\'e \ $k$ \ en (a,b), unitaire en \ $a$ \ et o\`u \ $R \in \A$ \ est un polyn\^ome en \ $a$ \ de degr\'e au plus \ $k-1$, \`a coefficients dans \ $\mathbb{C}[[b]]$, de valuation en \ $(a,b)$ \ sup\'erieur ou \'egale \`a \ $k-1$.\\
Dans ces conditions, l'\'el\'ement \ $P + b^2.R$ \ est unique et ne d\'epend que de l'id\'eal \ $I$ \ de \ $\A$. \\
De plus, on a l'\'egalit\'e dans \ $\A$
\begin{equation*}
 (-b)^k.B_E(-b^{-1}.a) = P  \tag{@}
 \end{equation*}
o\`u \ $B_E$ \ d\'esigne le polyn\^ome de Bernstein du (a,b)-module r\'egulier \ $E : = \A\big/I$. Ceci montre que \ $P \in \A$ \ ne d\'epend que de la classe d'isomorphisme de \ $E$.
\end{thm}

\parag{Remarque} Un (a,b)-module \ $E$ \ est un \ $\A-$module \`a gauche. Mais si l'on suppose que \ $E$ \ est r\'egulier (ou m\^eme seulement local, c'est \`a dire qu'il existe \ $N$ \ tel que \ $a^N.E \subset b.E$), il est complet pour la filtration \ $a-$adique et c'est naturellement un module \`a gauche sur l'alg\`ebre
$$ \hat{\mathcal{A}} : = \{ \sum_{n=0}^{\infty} \ S_n(a).b^n, S_n \in \mathbb{C}[[a]] \} $$
qui est la completion \ $a-$adique de \ $\A$. \\
Un isomorphisme de (a,b)-modules r\'eguliers est un isomorphisme de \ $\A-$modules \`a gauche et c'est alors  n\'ecessairement  un isomorphisme de \ $\hat{\mathcal{A}}-$modules \`a gauche. $\hfill \square$

\bigskip

\parag{D\'emonstration} Nous allons commencer par  montrer par r\'ecurrence sur le rang \ $k \geq 1$ \ du (a,b)-module monog\`ene r\'egulier \ $E$ \ l'assertion suivante :
\begin{itemize}
\item  pour tout g\'en\'erateur \ $x$ \ de \ $E$ \ il existe des nombres complexes \ $\lambda_1, \cdots, \lambda_k$ \ et des \'el\'ements \ $T_1, \cdots, T_k$ \ de \ $\mathbb{C}[[b]]$ \ v\'erifiant \ $T_j(0) = 1 \quad \forall j \in [1,k]$ \ tels que l'on ait
$$ \big[(a- \lambda_1.b).T_1.(a-\lambda_2.b).T_2\cdots (a - \lambda_k.b).T_k\big].x = 0 $$
dans \ $E$.
\end{itemize}
Supposons l'assertion montr\'ee en rang \ $k-1$. Consid\'erons alors une suite exacte
\begin{equation*}
0 \to F \to E \overset{\pi}{\to} E_{\lambda} \to 0  \tag{1}
\end{equation*}
o\`u \ $E_{\lambda}$ \ d\'esigne le (a,b)-module de rang \ $1$ \ (r\'egulier monog\`ene)   \ $\A\big/\A.(a - \lambda.b)$.\\
Le th\'eor\`eme d'existence des suites de Jordan-H{\"o}lder pour les (a,b)-modules r\'eguliers (voir [B.93]) assure de l'existence d'une telle suite exacte. Alors le (a,b)-module \ $F$ \ est r\'egulier de rang \ $k-1$. Comme \ $x$ \ est un g\'en\'erateur de \ $E$, \ $\pi(x)$ \ est un g\'en\'erateur de \ $E_{\lambda}$ \ et peut donc \^etre \'ecrit \ $ \pi(x) = S_k.e_{\lambda}$ \ o\`u \ $S_k \in \mathbb{C}[[b]]$ \ v\'erifie \ $S_k(0) = 1$ \ quitte \`a normaliser le g\'en\'erateur standard \ $e_{\lambda}$ \ de \ $E_{\lambda}$ \ v\'erifiant \ $(a - \lambda.b).e_{\lambda} = 0 $.\\
Pour faire avancer notre r\'ecurrence il nous suffit alors de montrer que l'\'el\'ement  \\ 
$y = (a - \lambda.b).S^{-1}_k.x$ \ est un g\'en\'erateur (sur \ $\A$) de \ $F$. D'abord c'est bien un \'el\'ement de \ $F = Ker\, \pi$ \ par d\'efinition de \ $S_k$. Alors il r\'esulte du lemme \ref{petit} que 
$y$ \ est  un g\'en\'erateur de \ $F$.\\
 Notre assertion est donc d\'emontr\'ee en posant \ $T_k : = S_k^{-1}$.\\
Notons \ $Q : = (a- \lambda_1.b).T_1.(a-\lambda_2.b).T_2\cdots (a - \lambda_k.b).T_k $ \ et consid\'erons le \ $\A-$module \ $ G : = \A\big/\A.Q$. On a une application \ $\A-$lin\'eaire \`a gauche surjective
$$ \varphi : G \to  E$$
d\'efinie en posant \ $\varphi(1) = x$ \ puisque, par construction, on a \ $Q.x = 0 $ \ dans \ $E$. Comme \ $T_1^{-1}\cdots T^{-1}_k.Q$ \ est un polyn\^ome unitaire en \ $a$ \ de degr\'e \ $k$ \ et de valuation \ $k$ \ en (a,b), \ $G$ \ est un \ $\mathbb{C}[[b]]-$module libre de rang \ $k$ \ (de base \ $1, a, \cdots, a^{k-1}$). L'application surjective \ $\varphi$ \ est donc un isomorphisme. On en conclut que tout id\'eal \`a gauche \ $I$ \  de \ $\A$ \ tel que \ $E : = \A\big/I$ \ soit un (a,b)-module r\'egulier est bien principal avec un g\'en\'erateur \ $x = P + b^2.R$ \ de la forme annonc\'ee.

\bigskip

Montrons la r\'eciproque. On suppose donc que \ $I = \A.x$ \ avec \ $x = P + b^2.R$ \ comme dans l'\'enonc\'e. Soit \ $y \in \A$ \ et \'ecrivons
$$ y = \sum_{n\geq 0} \ b^n.Y_n(a) $$
o\`u les \ $Y_n$ \ sont dans \ $\mathbb{C}[z]$. Posons alors \ $Y_n = Z_n.P + R_n$ \ pour chaque \ $n \geq 0$, o\`u \ $Z_n, R_n$ \ sont dans \ $\mathbb{C}[[b]][a]$ \ et \ $deg_a(R_n) \leq k-1 $. On aura alors
$$ y = \big(\sum_{n\geq 0}  b^n.Z_n\big).P + \sum_{n\geq 0} \ b^n.R_n .$$
Ceci montre que l'application \ $\mathbb{C}[[b]]-$lin\'eaire
$$ \psi : \mathbb{C}[[b]]^k \to  \A\big/\A.x  $$
d\'efinie par \ $\psi(S_0, \cdots, S_{k-1}) = [\sum_{j=0}^{k-1} \ S_j.a^j ]$ \ est surjective. \\
Montrons qu'elle est \'egalement injective. Supposons donc que
$$ \sum_{j=0}^{k-1} \ S_j.a^j = T.(P + b^2.R) .$$
Le cas \ $T = 0$ \ \'etant imm\'ediat, supposons que \ $T \not= 0 $. Alors on peut \'ecrire \\
 $T = b^N.T_1 + b^{N+1}.T_2$ \ o\`u \ $T_1 \in \mathbb{C}[a] \setminus \{0\}$ \ et \ $T_2 \in \A$.\\
 En raisonnant modulo \  $b^{N+1}.\A = \A.b^{N+1}$, on obtient que 
$$ \sum_{j=0}^{k-1} \ S_{j,N}.a^j = b^N.T_1(a).a^k  $$
o\`u \ $S_{j,N}$ \ d\'esigne la classe de \ $S_j$ \ dans \ $\mathbb{C}[[b]]\big/b^{N+1}.\mathbb{C}[[b]]$. On en d\'eduit que l'on a \ $T_1 = 0$. Contradiction. Donc  \ $\psi$ \ est injective.\\
Le quotient \ $\A\big/\A.x$ \ est donc  libre de rang \ $k$ \ sur \ $\mathbb{C}[[b]]$. C'est donc un (a,b)-module. Montrons qu'il est r\'egulier.\\
Pour cela, montrons que les \ $b^{-j}.a^j$ \ pour \ $j \in [0,k-1]$ \ engendrent (comme \ $\mathbb{C}[[b]]-$module) le satur\'e \ $E^{\sharp}$ \ de \ $E$ \ par \ $b^{-1}.a$ \ dans \ $E[b^{-1}]$. Il suffit, en fait, de montrer que l'on a 
 $$ b^{-k}.a^k.E \subset \sum_{j=0}^{k-1} \ b^{-j}.a^j.E .$$
 Ceci signifie que nous devons montrer l'inclusion
 $$ b^{-k}.a^k.\A \subset \sum_{j=0}^{k-1} \ b^{-j}.a^j.\A + b^{-k}.\A.x .$$
 Comme on a vu que \ $ \A = Im \psi + \A.x $ \ il suffit de voir que l'on a
 $$ a^k.\big(\sum_{j=0}^{k-1} \mathbb{C}[[b]].a^j \big) \subset  \sum_{j=0}^{k-1} \ b^{k-j}.a^j.\A + \A.x .$$
 Ceci va r\'esulter de l'assertion suivante :
   \begin{equation*} {\rm  Pour \  tout} \quad h \geq 0  \quad {\rm on \ a} \quad \quad  \mathbb{C}[[b]].a^{k+h} \subset \sum_{j=o}^{k-1} \ b^{k-j}.a^j. \mathbb{C}[[b]] + \A.x \tag{@@}
 \end{equation*}
 qui est clairement \'equivalente \`a 
  \begin{equation*} {\rm  Pour \  tout} \quad h \geq 0  \quad {\rm on \ a} \quad \quad  \mathbb{C}[[b]].a^{k+h} \subset \sum_{j=o}^{k-1}\ \mathbb{C}[[b]].b^{k-j}.a^j.  + \A.x \tag{@@bis}
 \end{equation*}
 Prouvons \ $(@@bis)$ \ par r\'ecurrence sur \ $h \geq 0$. Pour \ $h = 0$ \ on obtient
 $$ \mathbb{C}[[b]].a^k \subset \sum_{j=o}^{k-1}\ \mathbb{C}[[b]].b^{k-j}.a^j  + \A.x $$
 en \'ecrivant \ $x$ \ sous la forme
 $$ x = a^k + \sum_{j=0}^{k-1} \ S_j(b).b^{k-j}.a^j $$
 o\`u l'on a pos\'e
 $$ P : = a^k + \sum_{j=0}^{k-1} \ S_j(0).b^{k-j}.a^j \qquad {\rm et} \qquad R = \sum_{j=0}^{k-1} \ \frac{S_j(b) -S_j(0)}{b}.b^{k-j-1}.a^j .$$
 Supposons \ $(@@bis)$ \ montr\'ee pour \ $h \leq h_0$. Alors on a
 $$ \mathbb{C}[[b]].a^{k+h_0+1} \subset \ a.\mathbb{C}[[b]].a^{k+h_0} + \sum_{j=o}^{k-1}\ \mathbb{C}[[b]].b^{k-j}.a^j + \A.x $$
 en utilisant le fait que \ $\mathbb{C}[[b]].a \subset a.\mathbb{C}[[b]] + \mathbb{C}[[b]] $ \ et l'hypoth\`ese de r\'ecurrence. Celle-ci donne alors :
 $$ \mathbb{C}[[b]].a^{k+h_0+1} \subset \  a. \sum_{j=o}^{k-1}\ \mathbb{C}[[b]].b^{k-j}.a^j  +  \sum_{j=o}^{k-1}\ \mathbb{C}[[b]].b^{k-j}.a^j  + \A.x .$$
 On a
 $$a. \sum_{j=o}^{k-1}\ \mathbb{C}[[b]].b^{k-j}.a^j  \subset  \sum_{j=o}^{k-1}\ \mathbb{C}[[b]].b^{k-j}.a^j + \mathbb{C}[[b]].a^k \subset  \sum_{j=o}^{k-1}\ \mathbb{C}[[b]].b^{k-j}.a^j  + \A.x$$
 en utilisant le cas \ $h = 0 $ \ d\'emontr\'e plus haut. Ceci ach\`eve la preuve de la r\'egularit\'e de \ $E$.
 
 \bigskip
 
 Si maintenant on suppose que l'on a un autre g\'en\'erateur  \ $y = P_1 + b^2.R_1$ \ de \ $I$  \ tel que \ $P_1$ soit homog\`ene en (a,b) de degr\'e \ $k = rg(E)$, unitaire en \ $a$ \ et \  $R_1$ \ de valuation \ $\geq k-1$ \ en (a,b) et de degr\'e \ $\leq k-1$ \ en \ $a$, montrons que \ $P= P_1$ \ et \ $R = R_1$.\\
On a un \'el\'ement inversible    \ $u \in \A$ \ v\'erifiant  \ $ x = u.y$. Comme les inversibles de \ $\A$ \ sont de la forme \ $ \lambda + b.\xi$ \ o\`u \ $\lambda \in \mathbb{C}^*$ \ et \ $\xi \in \A$, on obtient, en supposant \ $u \not=1$ \  une \'egalit\'e
 $$ P + b^2.R = (1 + b^m.Q(a) + b^{m+1}.\eta).(P_1+ b^2.R_1) $$
 o\`u \ $\lambda = 1$ \ car \ $P$ \ et \ $P_1$ \ sont unitaires en \ $a$ \ de degr\'e \ $k$, et o\`u \ $m \geq 1$ \ est maximal \ tel que \ $u - 1 \in b^m.\A$. En regardant cette \'egalit\'e modulo \ $b^{m+1}.\A$ \ on obtient \ $P = P_1$ \ et \ $deg_a(b^m.Q(a).P_1) \leq k-1$ \ ce qui implique \ $Q(a) = 0$, contredisant le choix de \ $m$. On a donc \ $u = 1$ \ et \ $x = y$.
 
 \bigskip
 
  Montrons maintenant que les \ $(b^{-1}.a)^j, j \in [0,k-1]$, ou ce qui est \'equivalent, que les \ $b^{-j}.a^j, j \in [0,k-1]$, forment une base de \ $E^{\sharp}\big/b.E^{\sharp}$. Comme on sait d\'ej\`a que c'est un syst\`eme g\'en\'erateur, il suffit de prouver que c'est un syst\`eme libre. Mais comme on sait d\'ej\`a que le rang de \ $E$ \ est \'egal \`a \ $k$, le rang de \ $E^{\sharp}$ \ est aussi \'egal \`a \ $k$ \ et c'est n\'ecessairement une base. On en d\'eduit que le polyn\^ome caract\'eristique de \ $- b^{-1}.a$ \ agissant sur  \ $E^{\sharp}\big/b.E^{\sharp}$ \ est son polyn\^ome minimal. La relation \ $(@)$ \ r\'esulte alors de l'\'egalit\'e dans \ $E[b^{-1}]$ 
 $$ 0 = b^{-k}.x = b^{-k}.P + b^{-k+2}.R  $$
 en tenant compte du fait que \ $P$ \ est homog\`ene de degr\'e \ $k$ \ en $(a,b)$, unitaire en \ $a$ \ et que la valuation en \ $(a,b)$ \ de \ $R$ \ est au moins \'egale \`a \ $k-1$, alors que son degr\'e en \ $a$ \ est au plus \ $k-1$. $\hfill \blacksquare$

\bigskip

\subsection{\'El\'ement de Bernstein et  autres invariants.}

\begin{defn}[\'El\'ement de Bernstein]\label{Bernst.}
Soit \ $E$ \ un (a,b)-module monog\`ene r\'egulier de rang \ $k$, nous appellerons {\bf \'el\'ement de Bernstein de \ $E$}, que nous noterons \ $P_E \in \A$, l'\'el\'ement homog\`ene en (a,b) de degr\'e \ $k$, unitaire en \ $a$, d\'efini dans le th\'eor\`eme \ref{Structure}. Pour \ $I$ \ id\'eal \`a gauche de \ $\A$ \ tel que \ $E : = \A\big/I$ \ soit un (a,b)-module r\'egulier, l'\'el\'ement \ $x_I : = P_E + b^2.R_I \in \A$ \ sera appel\'e {\bf l'\'el\'ement caract\'erisitique} de \ $I$. $\hfill \square$
\end{defn}

\bigskip

Le lien entre les suites de Jordan-Hold{\"e}r de \ $E$ \ r\'egulier monog\`ene et l'\'el\'ement de Bernstein de \ $E$ \ est donn\'e par le corollaire suivant.

\begin{cor}\label{JH et P}
Soit \ $E$ \ un (a,b)-module monog\`ene r\'egulier de rang \ $k$ \ et soit
$$ 0 = F_0 \subset F_1\subset \cdots \subset F_k = E $$
une suite de Jordan-H{\"o}lder de \ $E$. Posons \ $F_j\big/F_{j-1} \simeq E_{\lambda_j}$ \ pour \ $j \in [1,k]$. \\
Alors on a \ $P_E = (a - \lambda_1.b)(a - \lambda_2)\cdots (a - \lambda_k.b)$. \\
Consid\'erons une autre suite de Jordan-H{\"o}lder \ $F'_j$ \ de \ $E$ \  et posons \ $F'_j\big/F'_{j-1} \simeq E_{m_j}$. Alors \ $(\mu_1, \cdots, \mu_k) \in \mathbb{C}^k$ \ sera dans l'orbite tordue de \ $(\lambda_1, \cdots, \lambda_k)$ \ sous l'action de \ $\mathfrak{S}_k$.
\end{cor}

\parag{Remarque} On notera que dans la situation du corollaire ci-dessus le polyn\^ome de Bernstein de \ $E$ \ est donn\'e par
$$ B_E(x) = (x + \lambda_1+k-1)(x + \lambda_2+k-2) \cdots (x + \lambda_k) .$$
En particulier \ $- \lambda_k$ \ est une racine du polyn\^ome de Bernstein de \ $E$. $\hfill \square$

\bigskip

\begin{lemma}[Compl\'ement] Soit  \ $E \simeq \A\big/I$ \ un (a,b)-module monog\`ene r\'egulier. Notons \ $x_I : = P_E + b^2.R_I$ \ l'\'el\'ement caract\'eristique de \ $I$. La partie homog\`ene \ $b^2.R_{k-1}$ \  en (a,b) de degr\'e \ $k+1$ \ de \ $b^2.R_I$ \ ne d\'epend, modulo \ $\mathbb{C}.[P_E,b] + \mathbb{C}.[P_E,a]$, que de la classe d'isomorphisme du (a,b)-module \ $E$.
\end{lemma}

\bigskip

On notera que les commutateurs \ $[P,a],[P,b]$ \ sont bien homog\`enes de degr\'e \ $k+1$ \ en (a,b) et de la forme \ $b^2.T$ \ o\`u \ $T$ \ est homog\`ene de degr\'e \ $k-1$ \ en (a,b).

\parag{Preuve}Gr\^ace \`a la proposition \ref{Iso} ceci revient \`a montrer, en posant \ $P : = P_E$, qu'une \'egalit\'e
\begin{equation*}
u.(P + b^2.R) = (P + b^2.S).v \tag{E}
\end{equation*}
o\`u \ $u,v$ \ sont des \'el\'ements de \ $\A$ \ qui sont \'egaux \`a \ $1$ \ modulo \ $\A.b + \A.a$ \ et \ $R,S$ \ des polyn\^omes en \ $a$ de degr\'e \ $\leq k-1$ \ et de valuation en (a,b) \ au moins \'egale \`a \ $k-1$, implique   \ $ b^2.R_{k-1} -  b^2.S_{k-1} \in \mathbb{C}.[P,a]+ \mathbb{C}.[P,b]$ \ ou \ $R_{k-1}$ \ et \ $S_{k-1}$ \ d\'esignent les parties homog\`enes en (a,b) de degr\'e \ $k-1$ \ de \ $R$ \ et \ $S$ \ respectivement.\\
Comme \ $P$ \ est  unitaire en \ $a$ \ on peut supposer que \ $ u = 1 + \lambda.b + {\lambda}'.a +  \xi $ \ et \ $v = 1 + \mu.b + {\mu}'.a +  \eta$ \ avec \ $\lambda, \mu \in \mathbb{C}$ \ et \ $\xi, \eta$ \ dans l'id\'eal bilat\`ere\ $  b^2\A + b.a\A + \A.a^2$. Alors la partie homog\`ene de degr\'e \ $k+1$ \ de l'\'egalit\'e \ $(E)$ \ donne
$$ \lambda.b.P + {\lambda}'.a.P + b^2.R_{k-1} = \mu.P.b + {\mu}'.a.P +  b^2.S_{k-1} ;$$
ceci implique \ ${\lambda}' = {\mu}'$, puis \ $\lambda = \mu$  \ et donc \ $b^2.R_{k-1} - b^2.S_{k-1} \in \mathbb{C}.[P,a] + \mathbb{C}.[P,b]$.$\hfill \blacksquare$

\bigskip

\begin{cor}\label{Invariants}
Tout (a,b)-module monog\`ene r\'egulier \ $E$ \ de rang \ $k$ \ est isomorphe \`a un quotient \ $\A\big/\A.x$ \ o\`u l'\'el\'ement \ $x$ \ de \ $\A$ \ est de la forme 
$$ x = P_E + b^3.S $$
avec \ $S$ \ un polyn\^ome en \ $a$ \ de degr\'e \ $\leq k-1$ \ et de valuation \ $\geq k-2$ \ en (a,b). Dans le choix d'un tel \ $x$, la partie homog\`ene \ $S_{k-2}$ \ de degr\'e \ $k-2$ \ de \ $S$ \ ne d\'epend que de la classe d'isomorphisme de \ $E$ \ modulo \ $\mathbb{C}.[P_E,(\alpha.a + b]$, o\`u \ $\alpha \in \mathbb{C}$ \ est tel que \ $deg_a([P_E,(\alpha.a + b)] \leq k-2$.
\end{cor}

\parag{Preuve} Le fait que \ $S_{k-2}$ \ ne d\'epende que de la classe d'isomorphisme de \ $E$ \ r\'esulte du lemme pr\'ec\'edent, puisque l'on a \ $b.P_E - P_E.b = b^2.k.a^{k-1} + Q$ \ avec \ $deg_a(Q)  \leq k-2 $. Il nous reste \`a montrer que tout \ $E$ \ est bien isomorphe \`a un tel quotient. Soit d\'ej\`a \ $y = P_E + b^2.R$ \ tel que \ $E \simeq \A\big/\A.y$. Soit \ $\mu \in \mathbb{C}$ \ le coefficient de \ $a^{k-1}$ \ dans la partie homog\`ene de degr\'e \ $k-1$ \ de \ $R$, et posons \ $\lambda = -\frac{\mu}{k}$. D\'efinissons \ $x = (1+\lambda.b).y.(1+ \lambda.b)^{-1}$. Alors on a un isomorphisme \ $\A-$lin\'eaire \`a gauche :
$$ \varphi : \A\big/\A.y \to \A\big/\A.x  $$
d\'efini par \ $\varphi(z) = z.(1 +\lambda.b)^{-1}$. Comme les id\'eaux \ $\A.x$ \ et \ $\A.y.(1+ \lambda.b)^{-1}$ \ sont \'egaux, \ $\varphi$ \ est bien un isomorphisme. Mais on a
$$ x = (1+\lambda.b).(P_E + b^2.R).(1 + \lambda.b)^{-1} $$
et on constate facilement que le coefficient de \ $b^2.a^{k-1}$ \ vaut \ $k.\lambda + \mu = 0$, ce qui prouve que \ $x = P_E + b^3.S$ \ avec \ $S$ \ un polyn\^ome de degr\'e \ $\leq k-1$ \ en \ $a$ \ et de valuation \ $\geq k-2$ \ en (a,b). $\hfill \blacksquare$

\parag{Exemples} 
\begin{enumerate}
\item Pour \ $P = a^2 - \alpha.ab + \beta.b^2$ \ on v\'erifie facilement que \ $[P,a]$ \ et \ $[P.b]$ \ engendrent l'espace vectoriel \ $\mathbb{C}.a.b^2 + \mathbb{C}.b^3$ \ si on a \ $4\beta \not= \alpha.(\alpha+2)$. Dans ce cas l'invariant \ $S_0$ \ du corollaire pr\'ec\'edent dispara\^it.

\item Si on consid\`ere, pour \ $\beta \in \mathbb{C}$ \ l'id\'eal \`a gauche de \ $\A$ \ engendr\'e par \ $a^2 + \beta.b^3$, on trouve la famille de (a,b)-modules de rang 2 qui sont not\'es \ $E_{1,0}(\alpha)$ \ dans la classification de la proposition 2.4 de [B.93] (avec \ $ \lambda = n = 1$). Montrons que l'invariant donn\'e par le corollaire ci-dessus (\'egal \`a \ $\beta$)  co{\"i}ncide avec l'invariant \ $\alpha$ \ de la classification :\\
Par d\'efinition \ $E_{1,0}(\alpha)$ \ admet une \ $\mathbb{C}[[b]]-$base \ $e_1, e_2$ \ dans laquelle on a 
$$ a.e_1 = (1 + \alpha.b).e_2 \quad {\rm et} \quad   a.e_2 = 0 .$$
Donc l'\'el\'ement \ $e_1$ \ engendre ce \ $\A-$module monog\`ene (r\'egulier) et il est annul\'e par
$$ x : = a.(1 + \alpha.b)^{-1}.a .$$
En fait il est plus pratique de consid\'erer \ $ (1 + \alpha.b)^{-1}.e_1$ \ comme g\'en\'erateur, car il est annul\'e par \ $y : = x.(1+ \alpha.b) =  a.(1 + \alpha.b)^{-1}.a.(1 + \alpha.b) $. Comme on a
$$ (1 + \alpha.b)^{-1}.a.(1 + \alpha.b) = a + \alpha.b^2 + b^3.\varphi(b) $$
o\`u \ $\varphi \in \mathbb{C}[[b]]$, l'annulateur du g\'en\'erateur \ $(1 + \alpha.b)^{-1}.e_1$ \ est l'id\'eal engendr\'e par l'\'el\'ement 
$$ y = a^2 + \alpha.ab^2 + ab^3.\varphi(b) .$$
Comme \ $P_E = a^2$ \ et que \ $[a^2,a] = 0, [a^2, b] = 2a.b^2 - 2b^3 $, l'invariant \ $S_0$ \ donn\'e par le corollaire est le coefficient de \ $b^3$ \ dans \ $\alpha.a.b^2 - \frac{\alpha}{2}.[a^2, b] = \alpha.b^3 $. Donc le corollaire et la classification des (a,b)-modules r\'eguliers de rang 2 nous donnent l'\'egalit\'e annonc\'ee.\\
A titre d'exercice sur les \'equations diff\'erentielles, le lecteur pourra d\'eterminer, dans la pr\'esentation du (a,b)-module \ $\A\big/\A.(a^2 + \alpha.b^3)$ \ sous la forme d'une \ $\mathbb{C}[[b]]-$base \ $e_1, e_2$ \ v\'erifiant \ $a.e_1 = e_2, \quad a.e_2 = -\alpha.b^3.e_1$ \ des \'el\'ements \ $S, T \in \mathbb{C}[[b]], T(0) = 1$ \ tels que l'\'el\'ement 
$$ \varepsilon : = S.e_1 + T.e_2 $$
satisfasse \ $a.\varepsilon = 0 $ (on sait \`a priori qu'il existe et qu'il est unique). $ \hfill \square$
\end{enumerate}

\bigskip

\subsection{Le second th\'eor\`eme de structure.}

\bigskip

\begin{thm}\label{Factorisation 1}
Soit \ $E$ \ un (a,b)-module monog\`ene r\'egulier de rang \ $k$. Soit \ $0 \subset F_1 \subset \cdots F_k = E$ \ une suite de Jordan-H{\"o}lder de \ $E$ \ et posons \ $F_j\big/F_{j-1} \simeq E_{\lambda_j}$ \ pour \ $j \in [1,k]$. Alors il existe \ $S_1, \cdots, S_{k-1} \in \mathbb{C}[b]$ \ v\'erifiant \ $S_j(0) = 1\quad \forall j$, tels que \ $E$ \ soit isomorphe \`a \ $\A\big/\A.\Pi$ \ o\`u 
\begin{equation*}
\Pi = (a - \lambda_1.b).S_1^{-1}.(a - \lambda_2.b).S_2^{-1}\cdots (a - \lambda_{k-1}.b).S_{k-1}^{-1}.(a - \lambda_k.b) \tag{**}
\end{equation*}
Si l'on suppose que pour chaque \ $j \in [1,k-1]$ \ on a \ $\lambda_j \not\in \lambda_{j+1} - \mathbb{N}$ \ alors on peut choisir les polyn\^omes \ $S_j$ \ v\'erifiant de plus \ $deg(S_j) \leq k-j-1 $.\\
Dans le cas g\'en\'eral, posons \ $m_j : = \sup \{ \lambda_h - \lambda_j, h \in [j+1,k]\big/  \lambda_h - \lambda_j \in \mathbb{N}\}$ si cet ensemble est non vide, et \ $m_j : = -1$ \ sinon.\\
 Alors on peut choisir \ $S_j$ \ de degr\'e \ $\leq k-j+m_j$.
\end{thm}

La d\'emonstration du th\'eor\`eme utilisera le lemme suivant.

  \begin{lemma}\label{tech.}
    Soit \ $E$ \ un (a,b)-module et soit \ $\pi : E \to E_{\mu}$ \ un morphisme surjectif. Soient \ $\lambda_1, \cdots, \lambda_k $ \ des nombres complexes, et soient \ $T_1, \cdots, T_k$ \ des inversibles de \ $\mathbb{C}[[b]]$.\\
    Notons \ $j_1, \cdots, j_h$ \ les \'el\'ements de \ $[1,k]$ \ pour lesquels on a \ $\mu \in \lambda_{j_i} - \mathbb{N}$. Posons \ $\mu = \lambda_{j_1} - m_1= \cdots = \lambda_{j_h} - m_h $ \ et \ $m : =  \sup_{i=1}^h \{m_i \}$ \ et \ $m : = -1$ \ s'il n'existe aucun\ $j$ \ tel que \ $\mu \in \lambda - \mathbb{N}$. Alors l'image de l'application induite par \ $\pi$
    $$ \pi' : \Big((a - \lambda_1.b).T_1.(a - \lambda_2.b) \cdots T_{k-1}.(a - \lambda_k.b).T_k\Big)(E) \to E_{\mu}  $$ 
   contient  \ $b^{k + m+1}.E_{\mu}$.\\
   Dans le cas o\`u  \ $\mu \not\in \lambda_j - \mathbb{N} \quad \forall j \in [1,k]$ \ alors l'image de \ $\pi'$ \ est exactement \ $b^k.E_{\mu}$.
    \end{lemma}
    
    \parag{Preuve}  Commen{\c c}ons par montrer le r\'esultat dans le cas o\`u \ $E = E_{\mu}$ \ et \ $\pi = Id$. \\
    Remarquons qu'il suffit de d\'emontrer ce lemme pour  \ $T_1 = \cdots = T_k = 1$ \ puisque l'on a \ $T.b^j.E = b^j.E$ \ pour tout \ $j \geq 0$ \ et tout \ $T \in \mathbb{C}[[b]]$ \ inversible.\\
   On a, pour chaque entier  \ $j \geq 0$ \ la relation
     $$ (a -\lambda_1.b)\cdots (a - \lambda_k.b).b^j = b^j.(a - (\lambda_1-j).b)\cdots (a - (\lambda_k-j).b)$$
     ce qui donne quand \ $a.e_{\mu} = \mu.b.e_{\mu}$ :
     $$(a -\lambda_1.b)\cdots (a - \lambda_k.b).b^j.e_{\mu} = (\mu - \lambda_1+j)\cdots (\mu - \lambda_k + j).b^{j+k}.e_{\mu} .$$
     Donc si \ $j \geq m+1$, on aura \ $(a -\lambda_1.b)\cdots (a - \lambda_k.b)(E_{\mu}) \supset b^{j+k}.e_{\mu}$ \ et donc l'inclusion 
      $$b^{k+m+1}(E_{\mu}) \subset (a -\lambda_1.b)\cdots (a - \lambda_k.b)(E_{\mu}). $$
      Dans le cas o\`u \ $\mu \not\in \lambda_j - \mathbb{N} \quad \forall j \in [1,k]$ \ ce raisonnement donne imm\'ediatement \ $(a -\lambda_1.b)\cdots (a - \lambda_k.b)(E_{\mu}) = b^k.E_{\mu}$.
    
       \bigskip
   
     Dans le cas g\'en\'eral la surjectivit\'e de l'application \ $\pi'$ \  induite par \ $\pi$ \ est cons\'equence imm\'ediate  du carr\'e commutatif suivant, de la surjectivit\'e de \ $\pi$ \ et des deux fl\`eches verticales \\
  $$  \xymatrix{ E \ar[d]^u \ar[r>]^{\pi} & E_{\mu} \ar[d]^u \\
                           u(E) \ar[r]^{\pi'} & u(E_{\mu})} $$
   o\`u l'on a pos\'e \ $u =  (a - \lambda_1.b).T_1.(a - \lambda_2.b) \cdots T_{k-1}.(a - \lambda_k.b).T_k \in \A$. $\hfill \blacksquare$

\parag{D\'emonstration du th\'eor\`eme}  Notons  \ $\pi_j$ \ la projection \ $F_j \to F_j \big/F_{j-1} \simeq E_{\lambda_j}$, pour \  $j \in [1,k]$. Comme \ $\pi_k$ \ est surjective, on peut choisir un g\'en\'erateur \ $x_k$ \  de \ $E$ \ de sorte que \ $\pi_k(x_k)$ \ soit un g\'en\'erateur (sur \ $\A$ ) \ $e_{\lambda_k}$ \ de \ $E_{\lambda_k}$ \ v\'erifiant \ $(a - \lambda_k.b).e_{\lambda_k} = 0 $. Alors on aura  \ $z_{k-1} : = (a - \lambda_k.b).x_k \in F_{k-1}$ \ et qui sera un g\'en\'erateur (sur \ $\A$ ) \ de \ $F_{k-1}$ \ gr\^ace au lemme \ref{petit}. \\
 Montrons par r\'ecurrence descendante sur \ $j \in [1,k]$ \ que l'on peut trouver les polyn\^omes \ $1 = S_k, \cdots, S_j, \cdots, S_1$ \ et des \'el\'ements \ $ \xi_j \in F_j, \forall j \in [1,k-1]$ \ tels que les \'el\'ements  \ $x_j : = x_k - \sum_{h=j}^{k-1} \xi_h $ \ soient des g\'en\'erateurs de \ $E$ \ et que pour chaque \ $j \in [1,k-1]$ \ (avec la convention \ $m_k = 0$) l'\'el\'ement
 $$ z_j : = (a - \lambda_{j+1}.b).S_{j+1}^{-1}\cdots S_{k-1}^{-1}.(a - \lambda_k.b).x_{j+1} $$
 soit un g\'en\'erateur de \ $F_j$ \ et que l'on ait
 $$ deg(S_j) \leq k - j + m_j  .$$

  Comme \ $S_k = 1$ \ et \ $x_k$ \ ont \'et\'e construits, supposons que \ $S_k, \cdots, S_{j+1}$ \ ainsi que \ $x_k, \cdots, x_{j+1}$ \ construits et construisons \ $S_j$ \ et \ $\xi_j$. Alors \ $z_j$ \ est un g\'en\'erateur de \ $F_j$ \ par hypoth\`ese. Posons \ $\pi_j(z_j) = T.e_{\lambda_j}$ \ o\`u  \ $T \in \mathbb{C}[[b]]$ \ et o\`u l'on peut supposer \ $T(0) = 1$, quitte \`a normaliser convenablement le g\'en\'erateur standard \ $e_{\lambda_j}$ \ de \ $E_{\lambda_j}$ \ v\'erifiant \ $a.e_{\lambda_j} = \lambda_j.b.e_{\lambda_j}$, puisque \ $\pi_j$ \ est surjective et que \ $z_j$ \ est un g\'en\'erateur de \ $F_j$.\\
  Utilisons maintenant le lemme \ref{tech.} qui nous donne l'inclusion
  $$ b^{k-j+m_j+1}.E_{\lambda_j} \subset \pi_j\big[(a - \lambda_{j+1}.b).S_{j+1}^{-1}.(a -\lambda_{j+2}.b).S_{j+2}^{-1} \cdots S_{k-1}^{-1}.(a - \lambda_k.b)(F_j)\big] .$$
  Posons alors \ $ T = T_1 + b^{k-j+m_j+1}.T_2 $ \ o\`u \ $T_1$ \ est un polyn\^ome de degr\'e \ $\leq k-j+m_j $ \ v\'erifiant \ $T_1(0) = 1$. Posons \ $S_j : = T_1$. \\
  Soit \ $\xi_j \in F_j$ \ tel que 
   $$\pi_j((a - \lambda_{j+1}.b).S_{j+1}^{-1}.(a -\lambda_{j+2}.b).S_{j+2}^{-1}\cdots S_{k-1}^{-1}.(a - \lambda_k.b)(\xi_j)) =  b^{k-j+m_j+1}.T_2.e_{\lambda_j}.$$
   Alors on aura
  $$ \pi_j((a - \lambda_{j+1}.b).S_{j+1}^{-1}.(a -\lambda_{j+2}.b).S_{j+2}^{-1}\cdots S_{k-1}^{-1}.(a - \lambda_k.b)(x_{j+1} - \xi_j)) = T_1.e_{\lambda_j} $$
 et donc 
  $$  \pi_j(S_j^{-1}.(a - \lambda_{j+1}.b).S_{j+1}^{-1}.(a -\lambda_{j+2}.b).S_{j+2}^{-1}\cdots S_{k-1}^{-1}.(a - \lambda_k.b)(x_{j+1} - \xi_j)) = e_{\lambda_j}$$
  ce qui montre que l'\'el\'ement
  $$ z_{j-1} : = (a - \lambda_j.b).S_j^{-1}.(a - \lambda_{j+1}.b).S_{j+1}^{-1}.(a -\lambda_{j+2}.b) \cdots S_{k-1}^{-1}.(a - \lambda_k.b)(x_{j+1} - \xi_j) $$
  est dans \ $ Ker \, \pi_j = F_{j-1}$. Comme \ $F_j \subset F_{k-1} \subset a.E + b.E $, \ $x_j : = x_{j+1} - \xi_j$ \ est bien un g\'en\'erateur de \ $E$. On a donc prouv\'e le pas de r\'ecurrence.\\
   L'assertion du th\'eor\`eme en d\'ecoule facilement puisque le morphisme surjectif \\
    $\A\big/\A.\Pi \to E$ \ entre deux (a,b)-modules de m\^eme rang est n\'ecessairement un isomorphisme. $\hfill \blacksquare$
   
   \parag{Remarque} Si on part d'un (a,b)-module monog\`ene r\'egulier \ $E$ \ qui poss\`ede un g\'en\'erateur \ $e$ \ dont l'annulateur est l'id\'eal \ $\A.\Pi$ \ o\`u \ $\Pi \in \A$ \ est donn\'e par l'\'equation \ $(^{**})$, les \ $S_1,\cdots,S_{k-1}$ \ \'etant des inversibles de \ $\mathbb{C}[[b]]$, on peut modifier le choix du g\'en\'erateur de mani\`ere que son annulateur soit l'id\'eal \ $\A.\Pi'$, o\`u \ $\Pi'$ \ est encore donn\'e par une \'equation \ $(^{**})$ \ avec maintenant des \ $S_j$ \ dans \ $\mathbb{C}[b]$ \ v\'erifiant \ $S_j(0) = 1$ \ et \ $deg(S_j) \leq k-j+m_j$.\\
   Ceci r\'esulte imm\'ediatement de la d\'emonstration du th\'eor\`eme.\\
   On prendra garde que l'on ne change pas la suite ordonn\'ee \ $ (\lambda_1, \cdots, \lambda_k)$ \ dans cette op\'eration.\\
    En particulier la suite de Jordan-Hold{\"e}r consid\'er\'ee  \  $0 \subset F_1 \subset \cdots \subset F_k = E$ \ ne change pas.  Le changement \'eventuel de suite de Jordan-Hold{\"e}r sera examin\'e au paragraphe suivant (voir la proposition \ref{reordonner} et le corollaire qui la suit). $\hfill \square$
  
\bigskip

\begin{cor} \label{Factorisation 2}
Soit \ $E = \A\big/I$ \ un (a,b)-module monog\`ene r\'egulier de rang \ $k$. Il existe des \'el\'ements \ $\theta_1, \cdots, \theta_k \in \mathbb{C}[[b]] $ \ tels que l'on ait
$$x_I : =  P_E + b^2.R_I = (a - b.\theta_1(b))\cdots (a - b.\theta_k(b)) . $$
\end{cor}

\bigskip

On remarquera que dans la situation du corollaire ci-dessus on a 
 $$P_E = (a -\theta_1(0).b)\cdots (a - \theta_k(0).b) .$$

\parag{Preuve} Appliquons le th\'eor\`eme pr\'ec\'edent et posons \ $T_j : = S_j\cdots S_{k-1} $. Comme on a pour \ $S \in \mathbb{C}[[b]$ \ inversible et \ $\lambda \in \mathbb{C}$ \ l'identit\'e
$$ S.(a - \lambda.b).S^{-1} = a - b.\theta(b) $$
o\`u \ $\theta(b) = \lambda + b.S'.S^{-1} $, on aura
\begin{align*}
& T_1.\Pi = T_1.(a - \lambda_1.b).T_1^{-1}.T_2.(a - \lambda_2.b).T_2^{-1}.T_3 \cdots T_{k-1}.(a - \lambda_{k-1}.b).T^{-1}_{k-1}.(a - \lambda_k.b) \\
& \qquad  = (a - b.\theta_1(b))\cdots (a - b.\theta_k(b))
\end{align*}
ce qui permet de conclure. $\hfill \blacksquare$

\parag{Remarque} Quand on choisit un isomorphisme \ $E \simeq \A\big/\A.x$ \ avec \ $x = P_E + b^3.S$, la factorisation du corollaire pr\'ec\'edent se fera avec la condition
$$ \sum_{j=1}^k \ \theta_j'(0) = 0 .$$

\bigskip

\begin{prop}\label{Multiplicativite du B}
Soit \ $E$ \ un (a,b)-module monog\`ene r\'egulier et soit
$$ 0 \to F \to E \to G \to 0 $$
une suite exacte de (a,b)-modules. Alors \ $F$ \ et \ $G$ \ sont monog\`enes r\'eguliers et on a l'\'egalit\'e dans \ $\A$
$$ P_E = P_F.P_G .$$
On en d\'eduit que le polyn\^ome de Bernstein de \ $E$ \ est donn\'e par la formule :
  $$B_E(x) = B_F(x - rg_G).B_G(x) .$$
\end{prop}

\parag{Preuve} Notons \ $\pi : E \to G$ \ la surjection \ $\A-$lin\'eaire \`a gauche de la suite exacte de l'\'enonc\'e. Soit \ $e \in E$ \ un g\'en\'erateur de \ $E$. Alors \ $\pi(e)$ \ est un g\'en\'erateur de \ $G$ \ qui est donc monog\`ene r\'egulier\footnote{Pour l'assertion " $E$ \ r\'egulier implique \ $F$ \ et \ $G$ \ r\'eguliers", voir [B.93].}. Notons \ $Q$ \ l'\'el\'ement caract\'eristique de l'annulateur de \ $\pi(e)$ \ dans \ $G$. Alors la forme initiale en (a,b) de \ $Q$ \ est \ $P_G$.\\
Montrons que \ $Q.e$ \ est un g\'en\'erateur de \ $F$. En effet si \ $y \in F$ \ \'ecrivons \ $y = \xi.e$ \ o\`u \ $\xi \in \A$. On a \ $\pi(y) = \xi.\pi(e) = 0$ \ et donc \ $\xi = \eta.Q$. Ceci montre que \ $y =\eta.Q.e$ \ et donc que \ $Q.e$ \ engendre \ $F$. Soit \ $R$ \ l'\'el\'ement caract\'eristique de son annulateur. On a donc \ $P_F$ \ qui est la forme initiale en (a,b) de \ $R$. On a \ $R.Q.e = 0$ \ et si \ $P$ \ est l'\'element caract\'eristique de l'annulateur de \ $e$ \ dans \ $E$, on aura donc \ $ R.Q = u.P$ \ avec \ $u \in \A$. Mais les degr\'es des formes initiales de \ $P,Q,R$ \ sont respectivement \ $rg(E),rg(G),rg(F)$ \ avec \ $rg(E) = rg(F) + rg(G)$ \ et elles sont unitaires de ces degr\'es en \ $a$. Ceci montre que la forme initiale de \ $u$ \ vaut \ $1$ \ et que l'on a \ $P_E = P_F.P_G$.\\
On obtient alors
$$ b^{-rg(E)}.P_E = b^{-rg(G)}.[b^{-rg(F)}.P_F].b^{rg(G)}.b^{-rg(G)}.P_G $$
et on conclut gr\^ace \`a l'identit\'e  \ $ b^{-j}.b^{-1}.a.b^j = b^{-1}.a - j \quad \forall j \in \mathbb{Z}. \hfill \blacksquare$

On notera que quand \ $G$ \ est de rang \ $1$, on aura \ $G \simeq E_{\lambda}$ \ et \ $B_G(x) = x + \lambda$. On retrouve ainsi la remarque qui suit le corollaire \ref{JH et P}.

\subsection{Changement de suite de Jordan-H{\"o}lder.}

Commen{\c c}ons par un lemme qui permettra dans presque tous les cas de permuter deux facteurs cons\'ecutifs.

\begin{lemma}\label{Permut.}
Soient \ $\lambda, \mu$ \ deux nombres complexes distincts  et posons \\
 $\delta : = \lambda - \mu \not= 0  $.
Soit \ $S \in \mathbb{C}[[b]]$ \ v\'erifiant \ $S(0) = 1 $ \ et supposons que si \ $\delta \in \mathbb{N}^*$ \ on a \ $S_{\delta} = 0 $\footnote{On a not\'e \ $S_{\delta}$ \ le coefficient de \ $b^{\delta}$ \ dans la s\'erie \ $S \in \mathbb{C}[[b]]$.} ; soit \ $U \in \mathbb{C}[[b]]$ \ l'unique solution de l'\'equation diff\'erentielle \ $ b.U' = \delta.(U - S) $ \ v\'erifiant \ $U(0) = 1$ \ et ne pr\'esentant pas le terme \ $b^{\delta}$ \ dans son d\'eveloppement quand \ $\delta \in \mathbb{N}^*$, alors on a l'identit\'e suivante dans  l'alg\`ebre \ $\A$
\begin{equation*} (a - \mu.b).S^{-1}.(a - (\lambda-1).b) = U^{-1}.(a - \lambda.b)U.S^{-1}.U.(a - (\mu-1).b).U^{-1} . \tag{@}
\end{equation*}
\end{lemma}

\parag{Preuve} Remarquons d\'ej\`a que \ $U$ \ est donn\'e, en posant \ $S : = \sum_{h=0}^{+\infty} S_h.b^h$,  par
$$ U(b) = \sum_{h=0}^{+\infty} \ \frac{\delta}{\delta - h}.S_h.b^h . $$

 En multipliant \`a droite et \`a gauche par \ $U$ \ il s'agit de calculer
$$U.(a - \mu.b).S^{-1}.(a - (\lambda-1).b).U .$$
Comme on a \ $a.U = U.a + b^2.U' $ \ et \ $b.U' = \delta.(U - S)$ \ cela donne
\begin{align*}
& U.(a - \mu.b).S^{-1}.(a - (\lambda-1).b).U = \\
& \quad \quad = ((a-\mu.b).U - \delta.b.(U - S)).S^{-1}.(U.(a-(\lambda-1).b) + \delta.b.(U-S)) \\
& \quad \quad  = ((a-\lambda.b).U.S^{-1} + \delta.b)(U.(a -(\mu-1).b) - \delta.b.S) \\
& \quad \quad  = (a - \lambda.b)U.S^{-1}.U.(a - (\mu-1).b) - (a- \lambda.b).\delta.b.U + \\
& \quad \quad \quad \quad   \delta.b.U.(a - (\mu-1).b) - \delta^2.b^2.S \\
& \quad \quad  = (a - \lambda.b)U.S^{-1}.U.(a - (\mu-1).b) +  \\
& \quad \quad \quad  \quad  + \lambda.\delta.b^2.U - a.\delta.b.U + \delta.b.U.a - \delta.(\mu-1).b^2.U - \delta^2.b^2.S 
\end{align*}
et on conclut, en utilisant les \'egalit\'es 
\begin{align*}
& a.(b.U) - (b.U).a = b^2.(b.U)' \quad {\rm et} \\
& \delta.b^2.(b.U)' = \delta^2.b^2.U - \delta^2.b^2.S + \delta.b^2.U 
\end{align*}
que l'\'el\'ement  \ $ \xi : =  \lambda.\delta.b^2.U - a.\delta.b.U + \delta.b.U.a - \delta.(\mu-1).b^2.U - \delta^2.b^2.S $ \ de \ $\A$ \ est nul : en effet on obtient, en utilisant la relation \ $\lambda - \mu = \delta$ :
$$ \xi  =  \delta^2.b^2.U + \delta.b^2.U - ( \delta^2.b^2.U - \delta^2.b^2.S + \delta.b^2.U ) - \delta^2.b^2.S = 0  .\qquad \qquad \blacksquare $$

\parag{Remarques} 
\begin{enumerate}[i)]
\item Si \ $S$ \ est un polyn\^ome en \ $b$ \ de degr\'e \ $d$ \ on notera que \ $U$ \  est alors \'egalement un polyn\^ome en \ $b$ \ de degr\'e \ $d$.
\item On peut renverser le processus car l'\'equation diff\'erentielle \ $b.U' = \delta.(U - S)$ \ donne, pour \ $V : = U^{-1}$ \ et \ $\Sigma : = S.U^{-2}$, que \ $V$ \ est  solution de l'\'equation \\
 $b.V' = -\delta.(V - \Sigma)$ \ et v\'erifie \ $V(0) = 1$. On notera que le fait qu'il existe une solution \ $V \in \mathbb{C}[[b]]$ \ implique  que pour \ $\delta \in - \mathbb{N}^*$, on a n\'ecessairement \ $\Sigma_{-\delta} = 0 \  ! \hfill \square$
\end{enumerate}

\bigskip

\begin{prop}\label{reordonner}
Consid\'erons l'\'el\'ement
  $$P : = S_0^{-1}.(a - \lambda_1.b).S_1^{-1}.(a - \lambda_2).S_2^{-1}\cdots (a - \lambda_{k-1}.b).S_{k-1}^{-1}.(a - \lambda_k.b).S_k^{-1}$$
   de \ $ \A$ \ o\`u les \ $S_0, \cdots, S_k$ \ sont dans \ $\mathbb{C}[[b]]$ \ et v\'erifient \ $S_j(0) = 1\quad \forall j \in [0,k]$, les \ $\lambda_j, j \in [1,k]$ \ \'etant des nombres complexes arbitraires. Notons \ $(\mu_1, \cdots, \mu_k)$ \ un \'el\'ement de l'orbite tordue de \ $(\lambda_1,\cdots, \lambda_k)$ \ v\'erifiant 
    \begin{equation*}
    \Re(\mu_1)+1 \leq \Re(\mu_2)+2 \leq \cdots \leq \Re(\mu_k)+k . \tag{@}
    \end{equation*}
     Alors il existe des \'el\'ements \ $T_0, T_1, \cdots, T_k$ \ dans \ $\mathbb{C}[[b]]$ \ v\'erifiant \ $T_j(0) = 1$ \ pour tout \ $ j \in [0,k]$ \ et
$$ P = T_0^{-1}.(a - \mu_1.b).T_1^{-1}.(a - \mu_2).T_2^{-1}\cdots (a - \mu_{k-1}.b).T_{k-1}^{-1}.(a - \mu_k.b).T_k^{-1} .$$
\end{prop}

\parag{Remarque} La preuve de la proposition montrera qu'il existe toujours un \'el\'ement de l'orbite tordue de \ $(\lambda_1, \cdots, \lambda_k)$ \ v\'erifiant la condition \ $(@)$. Elle montrera m\^eme qu'un tel \'el\'ement est unique quand les nombres \ $\lambda_1, \cdots, \lambda_k$ \ sont r\'eels. $\hfill \square$

\parag{Preuve} D'apr\`es le lemme pr\'ec\'edent \ref{Permut.}, si on a un produit de la forme \\ 
$(a - \mu.b).S^{-1}.(a - (\lambda-1).b) $ \ avec \ $\Re(\mu) > \Re(\lambda - 1) + 1$, on aura \ $\Re(\delta) < 0 $ \ et donc on peut le remplacer par \ $U^{-1}.(a - \lambda.b)U.S^{-1}.U.(a - (\mu-1).b).U^{-1} $ \ et maintenant on aura
$$ \Re(\lambda) \leq \Re(\mu-1) + 1 $$
ce qui signifie, si \ $i$ \ d\'esigne le rang de \ $(a - \mu.b)$ \ dans le produit donnant \ $P$,  qu' en faisant  agir la transposition \ $t_{i,i+1}$ \ par l'action tordue sur \ $(\lambda_1,\cdots, \lambda_k)$, on a r\'etabli l'ordre souhait\'e pour le produit entre les rangs \ $i$ \ et \ $i +1$.\\
Il est alors clair que l'on arrive ainsi \`a un \'el\'ement de l'orbite tordue qui satisfait la condition demand\'ee. Il reste alors \`a rejoindre l'\'el\'ement \ $(\mu_1, \cdots, \mu_k)$ \ donn\'e par une succession de transpositions de ce type. Mais d'apr\`es ce qui vient d'\^etre dit il reste seulement \`a consid\'erer des termes \ $(a - \mu.b).S^{-1}.(a - (\lambda-1).b) $ \ pour lesquels  \ $\Re(\mu) =  \Re(\lambda - 1) + 1$. Deux cas sont possibles : ou bien \ $\lambda = \mu$ \ et il n'y a rien \`a faire, ou bien \ $\lambda \not= \mu$ \ et alors \ $\delta \in i.\mathbb{R}$. Le lemme \ref{Permut.} permet \`a nouveau de conclure sans restriction dans ces cas. $\hfill \blacksquare$

\bigskip

Avant de donner un premier  corollaire de la proposition \ref{reordonner} rappelons que dans la classification des (a,b)-modules r\'eguliers de rang 2 donn\'ee dans [B. 93] proposition 2.4, appara\^it la famille de (a,b)-modules monog\`enes \ $E_{\lambda,\lambda-n}(\alpha)$ \ qui est par d\'efinition, pour \ $\lambda \in \mathbb{C}, \alpha \in \mathbb{C}^*$ \ et \ $n \in \mathbb{N}^*$
$$ E_{\lambda,\lambda-n}(\alpha) : = \A \big/ \A.(a - (\lambda-n).b)(1 + \alpha.b^n)^{-1}.(a - (\lambda-1).b)) .$$

\begin{lemma}\label{exception}
Soient \ $\lambda \in \mathbb{C}, \alpha \in \mathbb{C}^*$ \ et \ $n \in \mathbb{N}^*$. Alors le (a,b)-module  \ $E_{\lambda,\lambda-n}(\alpha)$ \ pour \ admet une {\bf unique} suite de Jordan-Hold{\"e}r 
$$ 0 \to E_{\lambda-n} \to E \to E_{\lambda-1} \to 0 $$
\end{lemma}

\parag{Preuve} Notons \ $e_1$ \ le g\'en\'erateur et \ $e_2 : = (1+ \alpha.b^n)^{-1}(a -(\lambda-1).b).e_1 $. On a ainsi une \ $\mathbb{C}[[b]]-$base de \ $E_{\lambda,\lambda-n}(\alpha)$ \ et l'op\'erateur \ $a$ \ est d\'efini par les formules :
$$ a.e_1 = (\lambda-1).b.e_1 + e_2 + \alpha.b^n.e_2, \quad \quad  a.e_2 = (\lambda-n).b.e_2 .$$
Nous allons d\'eterminer tous les sous-(a,b)-modules normaux de rang 1 de \ $E_{\lambda,\lambda-n}(\alpha)$. Cela revient \`a d\'eterminer les \'el\'ements \  $\varepsilon : = S(b).e_1 + T(b).e_2 $ \ v\'erifiant \ $a.\varepsilon = \mu.b.\varepsilon$ \ et \ $\varepsilon \not\in b.E_{\lambda,\lambda-n}(\alpha)$. On obtient imm\'ediatement les  relations
\begin{align*}
& b.S'(b) =  (\mu - \lambda +1).S(b)  \\
& S(b).(1 + \alpha.b^n) + b^2.T'(b) = (\mu-\lambda+n).b.T(b) \tag{*}
\end{align*}
On a donc \ $ p : = \mu - \lambda+1 \in \mathbb{N}$ \ et  \ $S(b) = c.b^p $ \ o\`u \ $c \in \mathbb{C}$.
Pour \ $p = 0$ \ l'\'equation \ $(^*)$ \ impose \ $c = 0$ \ et \ $T = 0 $, cas exclu. On a donc \ $p \geq 1$ \ et \ $S(0) = 0$, d'o\`u \ $T(0) \not= 0$. La r\'esolution de l'\'equation \ $(^*)$ \ avec \ $S(b) = c.b^p$ \ donne pour \ $p + n -1 = \mu - \lambda + n \not= 0$ \ que 
$$ T = c_1.b^{n+p-1} + \frac{c}{n}.b^{p-1} - c.\alpha.b^{n+p-1}.Log\, b $$
ce qui impose \ $c = 0$ \ puisque \ $T \in \mathbb{C}[[b]]$, et donc \ $S = 0 $ \ ce qui conduit \`a une contradiction puisque \ $\mu - \lambda + n \not= 0$.\\
Si l'on suppose \ $\mu = \lambda - n$ \ et donc \ $p = 1-n$, on retombe sur le cas \ $p = 0$ \ d\'ej\`a trait\'e. $\hfill \blacksquare$

\bigskip

\begin{cor}\label{Ordre J-H}
Soit \ $E$ \ un (a,b)-module monog\`ene r\'egulier tel que son \'el\'ement de Bernstein s'\'ecrive
$$ P_E : = (a - \lambda_1.b) \cdots (a - \lambda_k.b) .$$
Soit \ $(\mu_1, \cdots, \mu_k)$ \ un \'el\'ement de l'orbite tordue de \ $(\lambda_1, \cdots, \lambda_k)$ \ v\'erifiant la condition  \ $(@)$ \ de la proposition pr\'ec\'edente. Alors il existe un g\'en\'erateur \ $e$ \ de \ $E$ \ sur \ $\A$ \ dont l'annulateur est l'id\'eal \ $\A.Q$ \ avec 
 $$Q : = (a - \mu_1.b).T_1^{-1}.(a - \mu_2).T_2^{-1}\cdots (a - \mu_{k-1}.b).T_{k-1}^{-1}.(a - \mu_k.b) .$$
 Si, de plus,\ $E$ \ n'admet pas de sous-quotient de rang 2 isomorphe \`a \ $E_{\lambda,\lambda-n}(\alpha)$ \ avec \ $\lambda \in \mathbb{C}, \alpha \in \mathbb{C}^*$ \ et \ $n \geq 1$, alors pour tout \ $(\nu_1, \cdots, \nu_k)$ \ dans l'orbite tordue de \ $(\lambda_1, \cdots, \lambda_k)$ \ on peut trouver un g\'en\'erateur de \ $E$ \ sur \ $\A$ \ dont l'annulateur est \'egal \`a \ $\A.R$, avec 
 $$ R : = (a - \nu_1).X_1^{-1}.(a - \nu_2.b).X^{-1}_2 \cdots (a - \nu_{k-1}.b).X_{k-1}^{-1}.(a - \nu_k.b) $$
 o\`u les \ $X_j$ \ sont dans \ $\mathbb{C}[[b]]$ \ et v\'erifient \ $X_j(0) = 1$.
 \end{cor}
 
 \parag{Preuve} La premi\`ere assertion est une cons\'equence imm\'ediate de la proposition \ref{reordonner}. Pour montrer la seconde assertion, consid\'erons sous quelles circonstances nous ne pouvons pas \'echanger deux facteurs cons\'ecutifs, c'est \`a dire sous quelles conditions le lemme \ref{Permut.} ne peut s'appliquer \`a un produit 
  $$(a - \mu.b).S^{-1}.(a - (\lambda-1).b).$$
   On doit avoir \ $\delta = \lambda - \mu \in \mathbb{N}^* $ \ et \ $S_{\delta} \not= 0 $. On aura donc un sous-quotient ( \`a savoir \ $F_{j+2}\big/F_j \subset E\big/F_j$ ) qui est de rang 2, monog\`ene et dont le g\'en\'erateur est annul\'e par \ $(a - \mu.b).S^{-1}.(a - (\mu+\delta -1).b)$, avec la condition \ $S_{\delta} \not= 0$. Ceci correspond \`a une \ $\mathbb{C}[[b]]-$base \ $e_1, e_2$ \ v\'erifiant
 $$ a.e_1 = (\mu+\delta -1).b.e_1 + S(b).e_2 \quad \quad  a.e_2 = \mu.b.e_2 \quad {\rm avec} \quad S(0) = 1 .$$
 Soit \ $\delta.T(b) = \delta - 1 +  S(b) - S_{\delta}.b^{\delta}$, et soit \ $U \in \mathbb{C}[[b]]$ \ la solution de l'\'equation diff\'erentielle \ $ b.U' = \delta.(U - T) $ \ v\'erifiant \ $U(0) = 1$ \ et \ $U_{\delta} = 0$ \ (remarquer que \ $T_{\delta} = 0$). Posons alors
 $$ \varepsilon : = e_1 + \frac{U(b) - 1}{b}.e_2 .$$
 On a, 
 \begin{align*}
 &  a.\varepsilon =  (\mu+\delta -1).b.e_1 + S(b).e_2 + \mu.(U(b)-1).e_2 + b^2.\frac{b.U'(b) - U(b) +1}{b^2}.e_2 \\
 & \quad  \quad  =  (\mu+\delta -1).b.e_1 + (S(b) - \delta.T(b)+ \delta ).e_2 + (\mu + \delta -1).(U(b)-1).e_2 \\
 & \quad \quad  = (\mu+\delta -1).b.\varepsilon + (1 + S_{\delta}.b^{\delta}).e_2
  \end{align*}
  ce qui montre que ce (a,b)-module de rang 2 est isomorphe \`a \ $E_{\mu+\delta,\mu}(S_{\delta})$. Comme notre hypoth\`ese exclut cette possibilit\'e, le corollaire est d\'emontr\'e. $\hfill \blacksquare$

\bigskip
 
 \subsection{D\'etermination du polyn\^ome de Bernstein.}
 
 Le but de ce paragraphe est de determiner l'\'el\'ement de Bernstein d'un sous-(a,b)-module monog\`ene \ $F$ \ (necessairement r\'egulier) d'un (a,b)-module r\'egulier \`a partir de la connaissance d'un \'el\'ement non nul convenable de \ $\A$ \ annulant un g\'en\'erateur de \ $F$. Le r\'esultat principal que nous obtenons est le th\'eor\`eme suivant dont nous montrerons l'efficacit\'e sur des exemples concrets au paragraphe 5.

 \begin{thm}\label{Utile}
 Soit \ $E$ \ un (a,b)-module r\'egulier et soit \ $e \in E$ \ un \'el\'ement annul\'e dans 
 \ $E$ \ par un \'el\'ement \ $x \in \A$ \ dont la forme initiale en (a,b) \ est unitaire en \ $a$ \ et de degr\'e \ $k$. Supposons que \ $\A.e$ \ soit de rang au moins \ $k$ \  comme (a,b)-module. Alors le (a,b)-module monog\`ene r\'egulier \ $F : = \A.e$ \ est exactement de rang \ $k$ \ et il a pour \'el\'ement de Bernstein la forme initiale de \ $x$.
 \end{thm}
 
 \parag{Remarque} On prendra garde que, sous les hypoth\`eses du th\'eor\`eme, le \ $\A-$module \`a gauche \ $\A\big/\A.x$ \ n'est pas, en g\'en\'eral, un (a,b)-module r\'egulier. A priori seule son image \ $F$ \  dans \ $E$ \ par le morphisme \ $\A-$lin\'eaire \`a gauche \ $ \varphi $ \ d\'efini en posant \ $\varphi(1) = e$, est un (a,b)-module (monog\`ene) r\'egulier. $\hfill \square$
 
 \bigskip
 
 Commen{\c c}ons par deux  lemmes :
 
 \begin{lemma}[Division.]\label{Division}
Soit \ $P$ \ un \'el\'ement homog\`ene de degr\'e \ $k$ \ en (a,b), unitaire en \ $a$. Alors pour tout \ $X \in \A$ \ il existe \ $Q \in \A$ \ et \ $R \in \mathbb{C}[[b]][a]$ \ de degr\'e au plus \ $k-1$ \ en \ $a$ \ uniques v\'erifiant 
$$ X = Q.P + R .$$
\end{lemma}

\parag{D\'emonstration} L'unicit\'e est facile : si \ $Q.P + R = 0$ \ en raisonnant modulo \ $b.\A$ \ on obtient \ $Q_0.a^k + R_0 = 0$ \ o\`u \ $R_0 \in \mathbb{C}[a]$ \ est de degr\'e \ $\leq k-1$. On a donc \ $R_0 = 0 = Q_0$. Mais comme la multiplication par \ $b$ \ est injective dans \ $\A$, on se ram\`ene imm\'ediatement quand \ $(Q, R) \not=(0, 0)$ \ au cas o\`u \ $(Q_0, R_0) \not= (0, 0)$ \ ce qui contredit le raisonnement pr\'ec\'edent.\\
\'Ecrivons \ $X : = \sum_{\nu=0}^{\infty} X_m(a).b^m $ \ avec \ $X_m \in \mathbb{C}[a]$, et soit \ $\mu : \mathbb{N} \to \mathbb{N}$ \ une application croissante telle que l'on ait \ $deg(X_m) \leq \mu(m)$. \\
\'Ecrivons \'egalement 
$$ X : = \sum_{n=0}^{\infty} \ H_n  $$
o\`u \ $H_n$ \ est homog\`ene de degr\'e \ $n$ \ en \ $(a,b)$.\\
 La division euclidienne de \ $H_n$ \ par \ $P$ \ donne 
$$ H_n = Q_{n-k}.P + R_n $$
o\`u \ $R_n$ \ est homog\`ene de degr\'e \ $n$ \ en \ (a,b) et de degr\'e \ $\leq k-1$ \ en \ $a$ \ et \ $Q_{n-k}$ \ est homog\`ene en (a,b) de degr\'e \ $n-k$. On a donc \ $R_n \in b^{n-k+1}.\A$ \ pour \ $n \geq k$, et la s\'erie \ $ \sum_{n=0}^{\infty} \ R_n$ \ converge dans \ $\A$ \ vers un \'el\'ement \ $R$ \ qui est un polyn\^ome en \ $a$ \ de degr\'e \ $\leq k-1$ \ \`a coefficients dans \ $\mathbb{C}[[b]]$.\\
Montrons que la s\'erie \ $\sum_{n=k}^{\infty} \ Q_{n-k}$ \ converge \'egalement dans \ $\A$ \ vers un \'el\'ement que l'on notera \ $Q$. Pour cela il suffit de montrer que pour chaque \ $p  \in \mathbb{N}$ \ l'ensemble des \ $n \in \mathbb{N}$ \ tels que \ $Q_n \not\in b^p.\A$ \ est fini. Mais \ $Q_{n-k} \not\in b^p.\A$ \ implique \ $deg_a H_n \geq n-p+1$ \ et donc \ $\mu(p-1) \geq n-p+1$ puisqu'il doit exister un entier \ $q \leq p-1$ \ tel que \ $X_q$ \ soit de degr\'e \ $\geq n-p+1$. Ceci montre que l'on a \ $n \leq \mu(p-1) + p-1 $ \ dans ces conditions, ce qui prouve notre assertion.\\
On a alors clairement \ $ X = Q.P + R$ \ puisque la multiplication ( \`a droite par \ $P$) \ est continue dans \ $\A$. $\hfill \blacksquare$

\bigskip
 
  \begin{lemma}[Inversibilit\'e.]\label{Inversible}
 Soit \ $ X : = \sum_{n \geq 0} \ P_n(a).b^n $ \ un \'el\'ement de  \ $\tilde{\mathcal{A}}$ \ v\'erifiant \ $P_0(0) = 1$. Soit \ $E$ \ un (a,b)-module  r\'egulier de rang \ $k$. Alors l'action de \ $X$ \ sur \ $E$ \ est bijective et   pour \ $e \in E$, il existe un \'el\'ement \ $Y \in \tilde{\mathcal{A}}$, que l'on peut supposer polynomial en \ $a$ \ et de degr\'e en \ $a$ \ au plus \'egal \`a \ $k-1$, de terme constant \'egal \`a \ $1$, tel que l'on ait dans \ $E$ \ l'\'egalit\'e
 $$ Y.X.e = e .$$
  \end{lemma} 
 
 \parag{Preuve} La compl\'etion de \ $E$ \ pour la filtration \ $a-$adique, qui est cons\'equence de la r\'egularit\'e de \ $E$, montre que la suite d'endomorphismes de \ $E$ \ donn\'ee par l'action des  \ $T_N : = \sum_{n=0}^N \ (1-X)^n $ \ converge vers l'inverse de l'action de \ $X$ \ sur \ $E$. Celle-ci est donc bien bijective.\\ 
 Pour \ $e\in E$ \ consid\'erons le (a,b)-module monog\`ene  \ $F : =\tilde{\mathcal{A}}.e \subset E$. Il est r\'egulier de rang \ $l \leq k$,  et il existe un polyn\^ome \ $P$,  unitaire en \ $a$, de degr\'e \ $l$ \  \`a coefficients dans \ $\mathbb{C}[[b]]$ \ et de valuation \ $l$ \ en (a,b), tel que l'on ait un isomorphisme de  \ $\tilde{\mathcal{A}}-$modules\footnote{Ceci est d\'emontr\'e plus haut dans le th\'eor\`eme \ref{Structure}.}
 $$ \varphi :  \tilde{\mathcal{A}}\big/  \tilde{\mathcal{A}}.P \to F  $$
 d\'efini par \ $\varphi(1) = e$. On peut ainsi \'ecrire dans \ $F$, pour chaque \ $n \geq 0$,
 $$ a^{l+n}.e = \sum_{j=0}^{l-1}  S_{n,j}(b).a^j.e $$
 o\`u \ $S_{n,j} \in \mathbb{C}[[b]]$ \ est de valuation \ $\geq n$. En \'ecrivant alors
 $$ T_N.e = \sum_{j=0}^{l-1} \ T_{N,j}(b).a^j.e  $$
 on constate que pour chaque \ $j$ \ on \ $T_{N,j} - T_{N-1,j}$ \ qui est de valuation au moins \'egale \`a \ $N-l$, et donc la suite \ $ T_{N,j}$ \ converge dans \ $\mathbb{C}[[b]]$ \ vers un \'el\'ement \ $T_j$;  il reste \`a poser \ $Y : = \sum_{j=0}^{l-1} \ T_j(b).a^j $ \ pour achever la d\'emonstration. $\hfill \blacksquare$
 
 \parag{Remarque} On notera que le lemme ci-dessus est un r\'esultat plus pr\'ecis que la simple inversibilit\'e de \ $X$ \ dans \ $\hat{\mathcal{A}}$ \ le compl\'et\'e \ $a-$adique de \ $\A$, qui agit sur tout (a,b)-module r\'egulier (ou m\^eme local). $\hfill \square$
 
\bigskip

 La d\'emonstration utilisera \'egalement la proposition suivante, que l'on rapprochera du th\'eor\`eme \ref{Structure} :
 
 \begin{prop}
 Soit \ $ x = Q_k(a, b) + b.R_k(a, b) $ \ un \'el\'ement de l'alg\`ebre \ $ \A $, o\`u \ $Q_k$ \ est un polyn\^ome unitaire en $a$, homog\`ene de degr\'e $k \geq 1 $ \ en \ $(a, b)$, et o\`u l'\'el\'ement \ $R_k(a, b) \in  \A $ \ est de valuation en \ $(a, b)$ \ au moins \'egale \`a \ $k$.\\
 Alors le quotient \ $E : = \A \big/\A.x $ \ est un \ $\mathbb{C}[[b]]-$module libre de rang \ $k$. C'est un \ $(a, b)$-module r\'egulier.
 \end{prop}
 
 \parag{Preuve de la proposition} Le quotient \ $\A\big/\A.x $ \ est sans $b-$torsion. En effet,  si \ $b.y = z.x$, posons \ $z = Q(a) + b.t $ \ o\`u \ $Q \in \mathbb{C}[a]$ \ et o\`u \ $t \in \A$. On obtient alors \ $a^k.Q(a) = 0 $ \ dans \ $\mathbb{C}[a] \simeq \A\big/b.\A $. Et donc \ $Q = 0 $ \ et \ $z \in b.\A$. On en d\'eduit que \ $y \in \A.x$, d'o\`u notre assertion.\\
 Le fait que \ $1, a, \cdots, a^{k-1}$ \ forme une famille \ $\mathbb{C}[[b]]-$libre dans \ $E$ \ est imm\'ediat. Pour voir qu'elle est g\'en\'eratrice consid\'erons \ $z_0 \in \A$ \ et \'ecrivons \ $z_0 = u_0 + \zeta_0.x + b.z_1$ \ avec \ $u_0 \in \sum_{j=0}^{k-1} \mathbb{C}.a^j, \zeta_0 \in \mathbb{C}[a]$ \ et \ $z_1 \in \A$. Le fait qu'une telle \'ecriture existe r\'esulte imm\'ediatement de l'\'egalit\'e \ $\mathbb{C}[a] =  \sum_{j=0}^{k-1} \mathbb{C}.a^j + (a^k)$. \\
 On construit ainsi par r\'ecurrence sur \ $n \in \mathbb{N}$ \ des suites \ $(z_n)_{n\in \mathbb{N}}, (u_n)_{n\in \mathbb{N}}, (\zeta_n)_{n\in \mathbb{N}}$ \ respectivement dans \ $\A, \quad \sum_{j=0}^{k-1} \mathbb{C}.a^j $ \ et \ $\mathbb{C}[a]$ \ v\'erifiant
 $$ z_n = u_n + \zeta_n.x + b.z_{n+1} \quad \forall n \in \mathbb{N}.$$
 On voit alors facilement que l'on a \ $ z_0 = \sum_n b^n.u_n + (\sum_n b^n.\zeta_n).x $ \ o\`u par construction  \ $\sum_n b^n.u_n \in \sum_{j=0}^{k-1} \mathbb{C}[[b]].a^j $. \\ 
 
 Prouvons la r\'egularit\'e. L'\'egalit\'e
 $$ x = a^k + \sum_{j=1}^k \lambda_j.b^ja^{k-j} + b.R_k(a, b) $$
 donne dans \ $E$, pour \ $h \in [0, k-1]$
 $$ a^{k+h} + \sum_{j=1}^k \mu_{j,h}.b^ja^{k-j+h} \in b.\sum_{j=0}^k b^ja^{k-j}.E $$
 et donc, dans \ $E[b^{-1}]$
 $$ (b^{-1}a)^ka^h + \sum_{j=1}^k  \nu_{j,h} (b^{-1}a)^{k-j}a^h \in \sum_{j=0}^{k-1} (b^{-1}a)^j.E $$
 ce qui montre que
 $$ F : = \sum_{j=0}^{k-1} (b^{-1}a)^j.E $$
 est stable par \ $b$ \ et \ $b^{-1}a$. On en conclut que \ $E$ \ est r\'egulier. $\hfill \blacksquare $
 
 \parag{D\'emonstration du th\'eor\`eme \ref{Utile}} \'Ecrivons \ $x = in(P) + Q$ \ o\`u la valuation de \ $Q$ \ en (a,b) est au moins \ $k+1$\footnote{Attention, maintenant \ $Q$ \ n'est plus n\'ecessairement dans \ $b.\A$, on ne peut donc pas appliquer la proposition pr\'ec\'edente.}.  Comme \ $in(P)$ \ est unitaire de degr\'e \ $k$ \ en \ $a$, on peut \'ecrire
 $$ Q = X.in(P) + R  $$
 avec \ $R$ \ de degr\'e \ $\leq k-1$ \ en \ $a$ \ et \ $X, R$ \ de valuations repectives  en \ $(a,b)$ \ au moins \'egales \`a  \ $1$ \ et \ $k+1$. Ceci r\'esulte du lemme de division \ref{Division}.

 Ceci montre d\'ej\`a que l'on peut r\'e\'ecrire \ $R = b^{2}.S$ \ o\`u \ $S$ \ est de valuation en \ $(a,b)$ \ au moins \'egale \`a \ $k -1$.\\
 On a, dans \ $E$ \ 
 $$ (1 + X).in(P).e = - b^{2}.S.e $$
 et on peut appliquer le lemme \ref{Inversible} \`a l'\'el\'ement \ $in(P).e \in E$ \ et \`a l'\'el\'ement \ $(1+X)$ \ de \ $\tilde{\mathcal{A}}$, puisque \ $X$ \ est de valuation \ $ \geq1$ \ en \ $(a,b)$. On en d\'eduit l'existence de \ $Y \in \tilde{\mathcal{A}} $, polyn\^ome de degr\'e \ $\leq k-1$ \ en \ $a$, de terme constant \'egal \`a \ $1$, tel que l'on ait dans \ $E$
 $$in(P).e =  -Y.b^{2}.S.e \in b.E .$$
 On en d\'eduit que le rang de \ $E$ \ est au plus \'egal \`a \ $k$, donc \'egal \`a \ $k$.\\
 Soit \ $Z : = in(P) + Y.b^{2}.S$. C'est un polyn\^ome en \ $a$ \ qui annule \ $e$. Sa forme initiale en \ $(a,b)$ \ est \ $in(P)$ \  qui est de degr\'e \ $k$ \ et unitaire en \ $a$. Le quotient \ $\tilde{\mathcal{A}}\big/\tilde{\mathcal{A}}.Z$ \ est donc un (a,b)-module r\'egulier monog\`ene de rang \ $k$. De plus l'application \ $\tilde{\mathcal{A}}-$lin\'eaire \ $\varphi : \tilde{\mathcal{A}}\big/\tilde{\mathcal{A}}.Z \to E$ \ d\'efinie par \ $\varphi(1) = e$ \ est surjective. Son noyau est donc un (a,b)-module de rang nul. Donc \ $\varphi$ \ est un isomorphisme. On conclut alors gr\^ace au  th\'eor\`eme \ref{Structure} de la fa{\c c}on suivante :  soit \ $x_I : = P_E + b^2.R_I$ \ l'\'el\'ement caract\'eristique de l'id\'eal \ $I : = \A.Z$. Il existe un \'el\'ement unitaire\footnote{C'est-\`a-dire dont le terme constant vaut \ $1$.} \ $u \in \A$ \ tel que l'on ait \ $Z = u.x_I$. La comparaison des formes initiales en (a,b) donne alors \ $P_E = in(P)$ \ d'o\`u la formule cherch\'ee pour le polyn\^ome de Bernstein.
  $\hfill \blacksquare$
  
  \parag{Remarque} Dans la situation du th\'eor\`eme pr\'ec\'edent  on a toujours une surjection \ $\A-$lin\'eaire \`a gauche \ $\A\big/ \A.x \to E$ \ mais elle n'est pas,en g\'en\'eral, injective.
  
  \bigskip

\newpage

\section{D\'eveloppements asymptotiques standards et \\
 (a,b)-modules monog\`enes r\'eguliers g\'eom\'etriques.}

\subsection{\'Ecriture canonique.}

\begin{defn}
Un (a,b)-module r\'egulier \ $E$ \  est dit {\bf g\'eom\'etrique} si les valeurs propres de \ $b^{-1}.a$ \ agissant sur \ $E^{\sharp}\big/b.E^{\sharp}$ \ sont dans \ $\mathbb{Q}^{*+}$, o\`u \ $E^{\sharp}$ \ d\'esigne le satur\'e de \ $E$ \ par \ $b^{-1}.a$.
\end{defn}

On notera que le g\'en\'erateur \ $e_{\lambda}$ \ de \ $E_{\lambda}$ \ correspond au mon\^ome \ $s^{\lambda-1}$ \ et que \ $E_{\lambda}$ \ est g\'eom\'etrique si et seulement si la monodromie est unipotente et le mon\^ome est localement de carr\'e int\'egrable \`a l'origine.

\begin{defn}
Nous appellerons {\bf d\'eveloppement asymptotique standard} une serie formelle du type
$$ \varphi(s) = \sum_{j \in [0,n]} \sum_{\alpha \in A} \ T_{j.\alpha}(s).s^{\alpha}.\frac{(Log s)^j}{j!} $$
o\`u les \ $T_{\alpha,j}$ \ sont dans \ $\mathbb{C}[[s]]$, o\`u \ $n \in \mathbb{N}$ \ et o\`u \ $A$ \ est un sous-ensemble {\bf fini} de \ $\mathbb{Q} \,\cap\, ]-1,+\infty[ $.
\end{defn}

Nous noterons par \ $\Xi_N$ \ l'ensemble des DAS pour lesquels l'entier \ $n$ est dans \ $[0,N]$, et par \ $\Xi$ \ la r\'eunion des \ $\Xi_N$ \ quand \ $N$ \ d\'ecrit \ $\mathbb{N}$.

\smallskip

Nous consid\'ererons toujours les \ $\Xi_N$ \ ainsi que \ $\Xi$ \ comme des \ $\A-$modules \`a gauche, l'action de \ $a$ \ \'etant donn\'ee par multiplication par \ $s$, celle de \ $b$ par la primitive sans constante. Explicitement, si l'on pose \ $e_{\alpha,j} : = s^{\alpha}.\frac{(Log s)^j}{j!}  $ \ on aura
\begin{align*}
& a.e_{\alpha,j} = e_{\alpha+1,j} \\
& b.e_{\alpha,0} = \frac{1}{\alpha+1}.e_{\alpha+1,0}  \\
& b.e_{\alpha,j} = \frac{1}{\alpha+1}.e_{\alpha+1,j} - \frac{1}{\alpha+1}.b(e_{\alpha,j-1}) \quad \quad{\rm  pour} \quad j \geq 1
\end{align*}

  \bigskip
  
  \begin{lemma}\label{canon}
  Tout \ $\varphi \in \Xi_n$ \ admet une \'ecriture unique de la forme 
  \begin{equation*}
 \sum_{\alpha \in A_{\varphi}} \ S_{\alpha,j}(b).s^{\alpha}.\frac{(Log s)^j}{j!}  \tag{Can}
  \end{equation*}
  o\`u les \ $S_{\alpha,j}$ \ sont \ dans \ $\mathbb{C}[[b]]$, o\`u \ $A_{\varphi}$ \ est un sous-ensemble fini de \\
   $\mathbb{Q} \,\cap\, (]-1,+\infty[ \times [0,n])$, les conditions suivantes \'etant r\'ealis\'ees :
  \begin{enumerate}[i)]
\item L'ensemble \ $A_{\varphi}$ \ est satur\'e, c'est \`a dire que \ $(\alpha,j) \in A_{\varphi}$ \ implique \ $(\alpha,i) \in A_{\varphi}$ \ pour chaque \ $i \in [0,j]$.
\item Pour chaque \ $(\alpha,j) \in A_{\varphi}$ \ on a, ou bien \ $S_{\alpha,j} \equiv 0$, ou bien \ $S_{\alpha,j}(0) \not= 0 $.
\item Si \ $(\alpha,j) \in A_{\varphi} $ \ il existe \ $k \geq j$ \ tel que \ $S_{\alpha,k}(0) \not= 0$.
\item Si \ $(\alpha,j) \in A_{\varphi}$ \ v\'erifie \ $S_{\alpha,j} \not= 0$ \ et  si \ $(\alpha + p,j)\in A_{\varphi}$ \  avec \ $p \in \mathbb{N}^*$ \ on a  \ $S_{\alpha+p,j} \equiv 0$.
\end{enumerate}
\end{lemma}

\parag{Remarque} En combinant les conditions iii) et iv) ci-dessus, on obtient que si \ $(\alpha,j) \in A_{\varphi}$ \ v\'erifie \ $S_{\alpha,j} \not= 0$ \ et on a \ $(\alpha+p,j) \in A_{\varphi}$ \ avec \ $p \in \mathbb{N}^*$, alors il existe \ $k > j$ \ tel que \ $(\alpha+p,k) \in A_{\varphi}$ \ et \ $S_{\alpha+p,k} \not= 0 $. $\hfill \square$

\bigskip

\noindent L'\'ecriture donn\'ee dans le lemme pr\'ec\'edent pour un DAS \ $\varphi$ \ sera appel\'e {\bf l'\'ecriture canonique} de ce DAS.\\
 On dira que \ $\varphi$ \ est de {\bf poids \ $p$} \ si dans son \'ecriture canonique le cardinal de l'ensemble \ $A_{\varphi}$ \ est \'egal \`a \ $p$.
 
 \parag{Exemple} Pour \ $\alpha \in \mathbb{Q}^{*+}$ \   l'\'ecriture canonique de la fonction
 $$ \varphi = s^{\alpha}.Log s + s^{\alpha+1}.Log s $$
 est donn\'ee par \
 $ \varphi : = (1 +(\alpha+1).b).s^{\alpha}.Log s + \frac{1}{\alpha+1}. s^{\alpha+1}$, puisque l'on a \\ $(\alpha+1).b(s^{\alpha}Log  s) = s^{\alpha+1}.Log s - \frac{1}{\alpha+1}.s^{\alpha+1}.  \hfill \square$

\parag{Preuve} Comme pour \ $\alpha \in ]-1,0] \, \cap\, \mathbb{Q} $ \ les sous-$\A-$modules :
$$ \Xi_{\alpha,n} : = \sum_{j \in [0,n]} \sum_{p \in \mathbb{N}} S_{\alpha+p,j}(b).s^{\alpha+p}\frac{Log s)^j}{j!} $$
sont en somme directe dans \ $\Xi_n$, il suffit de prouver le lemme dans le cas o\`u \ $\varphi \in \Xi_{\alpha,n}$ \ o\`u \ $\alpha \in ]-1,0] \, \cap\, \mathbb{Q} $ \ sera fix\'e dans la suite.\\
Supposons que l'on ait \'ecrit 
$$ \varphi = \sum_{(\alpha+p,j) \in A} \ S_{\alpha+p,j}(b).s^{\alpha+p}.\frac{(Log s)^j}{j!} $$
les conditions i) \`a iv) \'etant remplies. Posons
$$ j_0 : = \sup\{ j \in [0,n] / \exists p \in \mathbb{N} \quad (\alpha+p,j) \in A \}.$$
Nous noterons \ $p_0$ \ le plus petit entier tel que \ $(\alpha+p_0, j_0) \in A$ \ et \ $S_{\alpha+p_0,j_0} \not= 0$. On notera qu'en fait les conditions iii) et  iv) impliquent qu'en fait \ $p_0$ \ est le seul entier \ $p$ \ tel que \ $(\alpha+p,j_0)$ \ soit dans \ $A$ \ avec \ $S_{\alpha+p,j_0} \not= 0$. C'est donc le plus grand entier \ $p$ \ tel que \ $(\alpha + p,j_0)$ \ soit dans \ $A$ \ d'apr\`es la condition iii).\\
Consid\'erons maintenant les \ $j < j_0$ \ pour lesquels on a \ $S_{\alpha+p_0,j} \equiv 0$. Ou bien ceci a lieu pour tous les \ $j \in [0,j_0-1]$, et dans ce cas nous poserons
$$ A' : = A \setminus \{(\alpha+p_0,0), \dots, (\alpha+p_0,j_0) \} $$
ou bien nous noterons par \ $j_1$ \ le plus grand entier dans \  $[0,j_0-1]$  \ tel que \ $S_{\alpha+p_0,j_1}(0) \not= 0 $, et nous poserons
$$ A' : = A \setminus \{(\alpha+p_0,j_1+1), \dots, (\alpha+p_0,j_0) \} .$$
Posons alors :
$$ \psi : = \varphi - S_{\alpha+p_0,j_0}(b).s^{\alpha_0+p_0}\frac{(Log s)^{j_0}}{j_0!} =  \sum_{(\alpha+p,j) \in A'} \ S_{\alpha+p,j}(b).s^{\alpha+p}.\frac{(Log s)^j}{j!}$$
et monrons que les conditions i) \`a iv) sont v\'erifi\'ees pour cette \'ecriture de \ $\psi$.\\
Les propri\'et\'es i) et ii) sont claires. Pour v\'erifier iii) il suffit de regarder le cas o\`u \ $(\alpha+p_0,j) \in A'$ \ v\'erifie \ $S_{\alpha+p_0,j} \equiv 0$, ce qui ne peut arriver que dans le cas o\`u \ $j_1$ \ existe. Mais alors on a \ $j \leq j_1$ \ et on peut prendre \ $k = j_1$.\\
La propri\'et\'e iv) est \'evidente car \ $A'$ \ est un sous-ensemble de \ $A$ \ et on n'a pas changer la valeur des \ $S_{\alpha,j}$.

\bigskip

Nous pouvons maintenant montrer par r\'ecurrence sur le poids de \ $A$, l'unicit\'e de l'\'ecriture canonique. Supposons donc l'unicit\'e montr\'ee pour une \'ecriture canonique de poids \ $\leq p-1$. On part donc de deux \'ecritures canoniques
\begin{align*}
& \varphi =  \sum_{(\alpha+p,j) \in A} \ S_{\alpha+p,j}(b).s^{\alpha+p}.\frac{(Log s)^j}{j!} \\
& \qquad =  \sum_{(\alpha+q,h) \in B} \ T_{\alpha+q,h}(b).s^{\alpha+q}.\frac{(Log s)^h}{h!} 
\end{align*}
Il est alors facile de voir que les entiers \ $j_0$ \ puis \ $p_0$ \ sont les m\^emes pour \ $A$ \ et \ $B$. En effet, \ $j_0$ \  correspond \`a l'exposant maximal de \ $Log s$ \ et \ $\alpha+p_0$ \ l'exposant minimal de \ $s$ \ pour  \ $(Log s)^{j_0}$. Il est aussi facile de v\'erifier  que pour construire \ $A'$ \ et \ $B'$ \ on enl\`eve le m\^eme ensemble \`a \ $A$ \ et \ $B$.\\
Il nous reste alors \`a montrer que 
\begin{align*}
&  \psi : = \varphi - S_{\alpha+p_0,j_0}(b).s^{\alpha_0+p_0}\frac{(Log s)^{j_0}}{j_0!}\\
& \psi_1 : = \varphi - T_{\alpha+p_0,j_0}(b).s^{\alpha_0+p_0}\frac{(Log s)^{j_0}}{j_0!}
\end{align*}
co{\"i}ncident pour pouvoir conclure gr\^ace \`a l'hypoth\`ese de r\'ecurrence. Comme \ $\psi$ \ et \ $\psi_1$ \ ne pr\'esentent plus que des termes dont les  puissances de \ $Log s$ \ sont au plus \'egale \`a \ $j_0-1$, l'\'egalit\'e cherch\'ee r\'esulte du fait que si \ $U \in \mathbb{C}[[b]]$ \ est telle que \ $U(b).s^{\alpha_0+p_0}\frac{(Log s)^{j_0}}{j_0!}$ \ ne pr\'esente plus que des termes  dont les puissances de \ $Log s$ \ sont  au plus \'egale \`a \ $j_0-1$, alors \ $U$ \ est nul. On a donc \ $S_{\alpha_0+p_0,j_0} = T_{\alpha_0+p_0,j_0}$ \ et \ $\psi = \psi_1$.\\
On a donc montr\'e l'unicit\'e de l'\'ecriture canonique.\\
L'existence d'une \'ecriture
$$ \varphi =  \sum_{(\alpha+p,j) \in A} \ S_{\alpha+p,j}(b).s^{\alpha+p}.\frac{(Log s)^j}{j!} $$
avec \ $A$ \ fini (mais arbitraire) est un exercice simple laiss\'e au lecteur. Pour passer \`a une \'ecriture dans laquelle \ $A$ \ v\'erifie les conditions i) \`a iv), faisons une r\'ecurrence sur la puissance maximale \ $j_0$ \ de \ $Log s$ \ qui apparait dans \ $\varphi$.\\
Consid\'erons donc les termes correspondant \`a la puissance maximale \ $j_0$ \  de \ $Log s$. Pour \ $S_{\alpha,j_0}(b) = b^p.T $ \ avec \ $T(0) \not= 0$, remplacons dans l'ensemble \ $A$ \ l'\'el\'ement \ $(\alpha,j_0)$ \ par \ $(\alpha +p,0), \cdots, (\alpha+p, j_0)$, et \'ecrivons
$$ T(b).b^p.s^{\alpha}.\frac{(Log s)^{j_0}}{j_0!}  = \sum_{h=0}^{j_0}  T_{\alpha+p,h}(b).s^{\alpha+p}.\frac{(Log s)^h}{h!}.$$
On notera que maintenant on a \ $T_{\alpha+p,j_0}(0) \not= 0$.\\
Par ailleurs on peut regrouper les diff\'erents termes de la forme \ $T_{\alpha+p,j_0}(b).s^{\alpha+p}.\frac{(Log s)^{j_0}}{j_0!} $ \ venant de diff\'erents \ $\alpha+p$ \ en constatant que l'on a, pour \ $p < q $
\begin{align*}
&  T_{\alpha+p,j_0}(b).s^{\alpha+p}.\frac{(Log s)^{j_0}}{j_0!} + T_{\alpha+q,j_0}(b)).s^{\alpha+q}.\frac{(Log s)^{j_0}}{j_0!} = \\
& \qquad  (T_{\alpha+p,j_0}(b) + b^{q-p}.T_{\alpha+q,j_0}(b)).s^{\alpha+p}.\frac{(Log s)^{j_0}}{j_0!}  + \sum_{h < j_0} termes (\frac{(Log s)^h}{h!}),  
\end{align*}
et \ $T_{\alpha+p,j_0}(b) + b^{q-p}.T_{\alpha+q,j_0}(b)$ \ est inversible dans \ $\mathbb{C}[[b]]$ \ d\`es que \ $T_{\alpha+p,j_0}(b) $ \ l'est.\\
On retranche \`a \ $\varphi$ \ le terme \ $T_{\alpha+p_0,j_0}(b).s^{\alpha+p}.\frac{(Log s)^{j_0}}{j_0!}$ \ obtenu apr\`es regroupement (donc \ $p_0$ \ est minimal), et on applique l'hypoth\`ese de r\'ecurrence sur \ $j_0$.\\
Il est alors facile de voir que l'on arrive \`a une \'ecriture canonique pour \ $\varphi$ \ en compl\'etant l'ensemble \ $A'$ \ donn\'e par l'hypoth\`ese de r\'ecurrence  pour qu'il contienne \ $(\alpha+p_0,0), \cdots,(\alpha+p_0,j_0)$.  $\hfill \blacksquare$

\bigskip

\subsection{Le th\'eor\`eme de r\'ealisation.}

Soit \ $V$ \ un espace vectoriel complexe de dimension finie. Nous munirons \ $\Xi_n\otimes_{\mathbb{C}} V $ \ de la structure de \ $\A-$module \`a gauche d\'efinie par \ $a.(\xi\otimes v) = (a.\xi)\otimes v$ \ et pour \ $S \in \mathbb{C}[[b]]$ \ par \ $ S(b).(\xi\otimes v) = (S(b).\xi)\otimes v$. On a alors le th\'eor\`eme de r\'ealisation suivant pour tout (a,b)-module g\'eom\'etrique.

\bigskip
  
\begin{thm} \label{Realisation}
Soit \ $\varphi \in \Xi_n$. Alors 
\begin{enumerate}
\item  \ $E : = \A.\varphi \subset \Xi_n$ \ est un (a,b)-module r\'egulier monog\`ene g\'eom\'etrique.
\item Il existe des { \bf polyn\^omes} \ $T_{\alpha,j} \in \mathbb{C}[b]$ \ v\'erifiant \ $T_{\alpha,j} (0) \not= 0 $ \ et tels que si \ $\psi : = \sum_{(\alpha,j) \in A} \ T_{\alpha,j}(b).s^{\alpha}.(Log\, s)^j $ \ le (a,b)-module \ $F : = \A.\psi \subset \Xi_n$ \  soit isomorphe \`a \ $E$.
\item Tout (a,b)-module  r\'egulier g\'eom\'etrique est isomorphe \`a un sous-(a,b)-module de \ $\Xi_n\otimes_{\mathbb{C}}  V$, pour un entier \ $n$ assez grand et \ $V$ \ un espace vectoriel complexe de dimension finie.
\end{enumerate}
\end{thm}

\parag{D\'emonstration} 
Montrons la propri\'et\'e 1 par r\'ecurrence sur le poids de \ $\varphi$. L'assertion \'etant triviale pour le poids \ $1$ \ supposons l'assertion montr\'ee pour les poids \ $\leq k-1$ \ et consid\'erons \ $\varphi$ \ de poids \ $k$. Soit \ $(\alpha,j) \in A_{\varphi}$ \ avec \ $j$ \ maximal l'exposant \ $\alpha$ \ \'etant  donn\'e. Remarquons que dans ces conditions on a \ $S_{\alpha,j}(0) \not= 0 $ \ d'apr\`es les propri\'etes ii) et iii) de l'\'ecriture canonique ; donc \ $S_{\alpha,j}(b)$ \ est un inversible de \ $\mathbb{C}[[b]]$. Consid\'erons alors \ $\psi : = (a - \alpha.b).S_{\alpha,j}^{-1}.\varphi$. Il est facile de voir que dans l'\'ecriture canonique de \ $\psi$ \ on a \ $A_{\psi} \subset A_{\varphi} + 1$ \ o\`u \ $A + 1 : = \{(\beta +1,j) / \ (\beta,j) \in A \}$, et que \ $(\alpha+1,j) \not\in A_{\psi}$. Donc le poids de \ $\psi$ \ est \ $\leq k-1$. Alors  l'hypoth\`ese de r\'ecurrence donne que \ $\A.\psi$ \ est (a,b)-module r\'egulier monog\`ene g\'eom\'etrique.\\
Mais la relation entre \ $\varphi$ \ et \ $\psi$ \ donne la suite exacte
$$ 0 \to \A.\psi \to \A.\varphi \overset{\pi}{\to} E_{\alpha+1}  \to  0 $$
o\`u \ $\pi(S_{\alpha,j}^{-1}.\varphi) : = e_{\alpha}$. On en d\'eduit que \ $\A.\varphi$ \ est r\'egulier g\'eom\'etrique puisque \ $ \alpha+1 > 0$ \ est rationnel.

\bigskip

L'assertion 2  est cons\'equence imm\'ediate du th\'eor\`eme \ref{Factorisation 1}.

\bigskip

Montrons l'assertion 3 par r\'ecurrence sur le rang de \ $E$. On se ram\`ene alors \`a montrer que si on a une suite exacte
$$ 0 \to F \to E \overset{\pi}{\to} E_{\alpha} \to 0  $$
avec \ $E$ \ g\'eom\'etrique (ce qui implique que  \ $\alpha $ \ est rationnel et strictement positif) et \ $F$ \ isomorphe \`a un sous-module de  \ $\Xi_n\otimes V$ \ via un plongement \ $\A-$lin\'eaire \ $g : F \to \Xi_n\otimes V$, gr\^ace \`a l'hypoth\`ese de r\'ecurrence, alors \ $E$ \ est isomorphe \`a un sous-module de  \ $\Xi_{n+1}\otimes W$. Soit \ $e \in E$ \ un \'el\'ement de \ $E$ \ v\'erifiant \ $\pi(e) = e_{\alpha}$. Alors \ $z : = (a - \alpha.b).e $ \ qui est dans \ $F$. Posons \ $g((a - \alpha.b).e) = \sum_{p=1}^q \psi_p\otimes v_p$. On voit donc que le probl\`eme consiste alors \`a trouver  des \ $\varphi_p \in \Xi_{n+1}, p \in [1,q]$ \ tels que \ $(a - \alpha.b).\varphi_p = \psi_p$. C'est un exercice simple de montrer que de tels \ $\varphi_p$ \ existent toujours. On d\'efinit alors \ $W : = V \oplus \mathbb{C}.\varepsilon$ \ et une application \ $\A-$lin\'eaire \ $f :  E \to \Xi_{n+1}\otimes (V \oplus \mathbb{C}.\varepsilon)$ \ en posant \ $ f(e) = \sum_{p=1}^q \varphi_p\otimes v_p + s^{\alpha-1}\otimes  \varepsilon$, et \ $f(x) = g(x) $ \ pour \ $x \in F$, sous r\'eserve de v\'erifier que si \ $u \in \A$ \ v\'erifie \ $u.e = y \in F$ \ on a bien \ $u.f(e) = g(y) $ \ dans \ $\Xi_{n+1}\otimes W$. Mais notre construction montre que l'id\'eal annulateur de \ $e$ \ dans \ $E/F \simeq E_{\alpha}$ \ est \ $\A.(a-\alpha.b)$. On a donc \ $u = u_1.(a - \alpha.b)$, et on a donc \ $u.e = u_1.z$ 
\begin{align*}
& u.f(e) = u_1.(a -\alpha.b)\big[\sum_{p=1}^q \varphi_p\otimes v_p + s^{\alpha-1}\otimes \varepsilon\big] \\
& \qquad = \sum_{p=1}^q (u_1.\psi_p)\otimes v_p = g(u_1.z)
\end{align*}
puisque \ $g(z)  = \sum_{p=1}^q \psi_p\otimes v_p $. $\hfill \blacksquare$

\bigskip

\parag{Remarque} On notera que si on omet le terme \ $s^{\alpha-1}\otimes \varepsilon$ \ dans la preuve pr\'ec\'edente, on obtient bien une application \ $\A-$lin\'eaire \ $f$ \  de \ $E$ \ dans \ $\Xi_{n+1}\otimes V$ \ dont la restriction \`a \ $F$ \ est bien \ $g$, mais il est possible que \ $f$ \ ne soit pas injective. C'est par exemple le cas pour \ $E : = E_{\lambda+k,\lambda} \simeq \mathbb{C}[[b]].e_1 \oplus \mathbb{C}[[b]].e_2 $ \ avec \ $k\in \mathbb{N}^*,  \lambda-1 \in \mathbb{Q}^{+*}$ \ o\`u \ $a$ \ est d\'efini par les relations 
$$ a.e_1 = e_2 + (\lambda+k-1).b.e_1, \qquad a.e_2 = \lambda.b.e_2$$
 et o\`u  \ $g(e_2) : = s^{\lambda-1} \in \Xi_0$. En effet la relation \ $(a - (\lambda+k-1).b).e_1 = e_2$ \ imposera de choisir \ $f(e_1) = c.s^{\lambda+k-2} - \frac{\lambda-1}{k}.s^{\lambda-2}$. On constate alors que l'image de \ $f$ \ est contenue dans \ $\mathbb{C}[[b]].s^{\lambda-2} \simeq E_{\lambda-1}$ \ qui est de rang 1.\\
 On notera que nous avons choisi un (a,b)-module {\bf monog\`ene} g\'eom\'etrique dans l'exemple pr\'ec\'edent, pour bien montrer que le probl\`eme ne vient pas seulement du cas o\`u l'on cherche \`a \'etendre \`a \ $E_{\lambda} \oplus E_{\lambda}$ \ le plongement standard de \ $E_{\lambda}$ \ dans \ $\Xi_0$. $\hfill \square$

\subsection{Exposants.}

\bigskip

\begin{defn}
Soit \ $E$ \ un (a,b)-module g\'eom\'etrique et soit \ $\Lambda \in \mathbb{Q}/\mathbb{Z}$. Nous appellerons $\Lambda-$exposant de \ $E$, not\'e \ $exp_{\Lambda}(E)$ \  le plus petit \ $\lambda \in \Lambda$ \ tel que \ $b^{-1}.a - \lambda$ \ ne soit pas injectif sur \ $E^{\sharp}\big/b.E^{\sharp}$ \ o\`u \ $E^{\sharp}$ \ d\'esigne le satur\'e de \ $E$ \ par \ $b^{-1}.a$.
\end{defn}

\bigskip

On notera que pour \ $\lambda \in \Lambda$ \  la condition \ $\lambda < exp_{\Lambda}(E)$ \ implique que \ $-\lambda$ \ n'est pas racine du polyn\^ome de Bernstein  \ $B_E$, puisque ce dernier est le polyn\^ome caract\'eristique de l'action de \ $-b^{-1}.a$ \ sur \ $E^{\sharp}\big/b.E^{\sharp}$.

\bigskip

\begin{lemma}
Soit \ $E = \A.e$ \ un (a,b)-module monog\`ene g\'eom\'etrique.\\
 Soit \ $\lambda \in  \mathbb{Q}^{*+}$ \ v\'erifiant \ $\lambda < exp(E)$ \ et soit \ $S$ \ un \'el\'ement inversible de \ $\mathbb{C}[[b]]$. Posons \ $F : = \A.(a -\lambda.b).S^{-1}.e \subset E$. Alors on a \ $rg(F) = rg(E)$.
\end{lemma}

\parag{Preuve} Quitte \`a remplacer \ $e$ \ par \ $S^{-1}.e$ \ on peut supposer que \ $S = 1$. Soit \ $Q$ \ un g\'en\'erateur "standard" (voir le th\'eor\`eme \ref{Structure}) de l'annulateur de \ $(a- \lambda.b).e$ \ dans \ $F$ \ dont la forme initiale est de degr\'e \ $l = rg(F)$ \ et est unitaire en \ $a$. On aura \ $Q.(a-\lambda.b).e = 0$ \ dans \ $E$, et donc \ $Q.(a-\lambda.b) \in I$, l'annulateur de \ $e$. Ceci implique que les \'el\'ements de Bernstein \ $ P_F$ \ et \ $P_E$ \ v\'erifient une \'egalit\'e de la forme
$$ P_F.(a-\lambda.b) = R.P_E $$
o\`u \ $R$ \ est homog\`ene de degr\'e \ $k - l$ \ et unitaire en \ $a$. Si on a \ $l < k$, alors on aura \ $R = 1$ \ (on notera qu'il est \'evident \`a priori, d'apr\`es le lemme \ref{petit} que le rang de \ $F$ \ est \'egal \`a \ $k$ \ ou \ $k-1$). Posons \ $P_E = (a - \lambda_1.b) \dots (a-\lambda_k.b)$. D'apr\`es la proposition  \ref{homogene} ceci montre qu'il existe \ $j \in [1,k]$ tel que \ $\lambda_j+k-j = \lambda$, contredisant notre hypoth\`ese. On a donc \ $R$ \ qui est de degr\'e \ $1$ \ et \ $l = k$ . $\hfill \blacksquare$

\bigskip

\begin{defn}\label{nilpotence}
Soit \ $\varphi \in \Xi$ \ et soit  \ $A_{\varphi}$ \ l'ensemble satur\'e intervenant dans l'\'ecriture canonique de \ $\varphi$ \ (voir le lemme \ref{canon}). Notons \ $\tilde{\pi} : A_{\varphi} \to \mathbb{Q}/\mathbb{Z} $ \ la compos\'ee de l'application de projection \ $\pi : A_{\varphi} \to \mathbb{Q}^{*+}$ \ et de l'application quotient. Pour \ $\Lambda \in \tilde{\pi}$ \ posons
$$ n(\Lambda) : = \sup\{ j \in \mathbb{N} \ \big/ \ \exists (\alpha,j) \in A_{\varphi} \quad {\rm avec} \quad \alpha \in \Lambda \quad {\rm et} \quad S_{\alpha,j}(0) \not= 0\}.$$
Posons \'egalement pour \ $\Lambda \in \tilde{\pi}(A_{\varphi})$
$$ \Lambda^+ : = \sup \{ \alpha\in \Lambda \  \big/ \  (\alpha, n(\Lambda)) \in A_{\varphi} \quad {\rm et} \quad S_{\alpha,n(\Lambda)}(0) \not= 0 \}.$$
\end{defn}

\bigskip

\begin{prop}\label{exposants}
Soit \ $\varphi \in \Xi$ \ et soit
$$ \varphi = \sum_{(\alpha,j) \in A} \ S_{\alpha,j}(b).s^{\alpha}.\frac{(Logs)^j}{j!} $$
son \'ecriture canonique. Alors le rang de \ $E : = \A.\varphi \subset \Xi$ \ est donn\'e par la formule suivante
$$ rg(E) = \sum_{\Lambda\in \tilde{\pi}(A_{\varphi})} \  n(\Lambda) .$$

 Pour chaque \ $\Lambda \in \tilde{\pi}(A_{\varphi})$ \ le polyn\^ome de Bernstein de \ $E$ \ contient exactement \ $n(\Lambda)$ \ racines dans \ $\Lambda$ \ en comptant les multiplicit\'es. De plus il admet n\'ecessairement \ $-\Lambda^+$ comme racine.
\end{prop}

On voit donc que l'on peut lire sur le polyn\^ome de Bernstein dans le cas d'un (a,b)-module monog\`ene g\'eom\'etrique une bonne partie des informations sur les mon\^omes intervenant dans l'\'ecriture canonique d'un g\'en\'erateur.

\parag{Preuve} Nous allons montrer les assertions de la proposition dans ce cas par r\'ecurrence sur \ $j_0$ \ la puissance maximale des logarithmes intervenant dans le d\'eveloppement de \ $\varphi$.\\
 Pour \ $j_0 = 0$ \ l'\'ecriture canonique de \ $\varphi$ \ se r\'eduit \`a
$$ \varphi = \sum_{\alpha \in \pi(A)} S_{\alpha,0}(b).s^{\alpha}  $$
o\`u l'on a \ $S_{\alpha,0}(0) \not= 0$ \ pour chaque \ $\alpha \in \pi(A)$, et si \ $\alpha \not= \beta$ \ sont dans \ $\pi(A)$, ils ne sont pas congrus modulo \ $\mathbb{Z}$ \ \`a cause des  propri\'et\'es iii) et iv). Le rang est alors clairement le cardinal de \ $A$ \ ce qui correspond bien \`a la formule de l'\'enonc\'e dans ce cas.\\
De m\^eme le calcul de l'\'el\'ement de Bernstein de \ $E : = \A.\varphi$ \ est imm\'ediat dans ce cas, et on en d\'eduit que \ $B_E(x) = \prod_{\alpha \in \pi(A)} (x+\alpha+1) $.\\

\smallskip

Supposons donc la proposition d\'emontr\'ee quand l'exposant maximal des logarithmes est \ $\leq j_0-1$, et montrons-la quand \ $\varphi$ \ a des exposants de logarithmes major\'es par \ $j_0$.\\
Consid\'erons d\'ej\`a le cas o\`u on a un seul  terme \ $S_{\alpha,j_0}(b).s^{\alpha}.\frac{(Logs)^{j_0}}{j_0!}$ \ de l'\'ecriture canonique qui a l'exposant \ $j_0$ \ pour le logarithme. Quitte \`a consid\'erer \ $S_{\alpha,j_0}^{-1}(b).\varphi$, ce qui ne change pas le (a,b)-module \ $E$ \ consid\'er\'e, on peut supposer que l'on a
$$ \varphi = s^{\alpha}.\frac{(Logs)^{j_0}}{j_0!} + \psi $$
o\`u \ $\psi$ \ ne pr\'esente plus de logarithme avec exposant \ $j_0$. Alors
$$ \varphi_1 : = (a - (\alpha+1).b)[\varphi] $$
satisfait l'hypoth\`ese de r\'ecurrence. Donc le (a,b)-module \ $F : = \A.\varphi_1$ \ a son rang donn\'e par la formule 
$$ rg(F) = \sum_{M \in \tilde{\pi}(A_{\varphi_1}) }  \ n(M) $$
o\`u \ $A_{\varphi_1}$ \ indexe l'\'ecriture canonique de \ $\varphi_1$ :
$$\varphi_1 = \sum_{(\beta,j)\in A_{\varphi_1}} \ T_{\beta,j}(b).s^{\beta}.\frac{(Logs)^{j}}{j!}.$$
 Mais on constate facilement que pour \ $\Lambda \not= [\alpha] $ \ on a \ $n(\Lambda)_{\varphi_1} = n(\Lambda)_{\varphi}$. Montrons que pour \ $\Lambda = [\alpha]$ \ on a \ $n(\Lambda)_{\varphi_1} = n(\Lambda)_{\varphi} - 1$. L'in\'egalit\'e \ $\leq$ \ est claire puisque \ $\varphi_1$ \ n'a plus que des exposant \ $\leq j_0-1$ \ en \ $Log s$. Par ailleurs 
 $$(a - (\alpha+1).b)(s^{\alpha}.\frac{(Logs)^{j_0}}{j_0!})= \frac{s^{\alpha+1}}{\alpha+1}.\frac{(Logs)^{j_0-1}}{(j_0-1)!} $$
 montre que \ $(\alpha+1,j_0-1)\in A_{\varphi_1}$, et que \ $T_{\alpha+1,\j_0-1}(0) \not= 0 $. On a donc dans ce cas \ $n(\Lambda)_{\varphi_1} = j_0-1$. \\
 Une r\'ecurrence facile sur le nombre de \ $\Lambda \in \tilde{\pi}(A_{\varphi})$ \ pour lesquelles on a \ $n(\Lambda) = j_0$ \ permet d'achever la d\'emonstration du calcul du rang.\\
 La preuve ci-dessus montre que pour \ $\Lambda \in \tilde{\pi}(A_{\varphi})$ \ le nombre \ $- \Lambda^+$ \ est racine du polyn\^ome de Bernstein de \ $E : \A.\varphi$. L'assertion sur le nombre de racines dans  \ $\Lambda \in \tilde{\pi}(A_{\varphi})$ \ d\'ecoule de la r\'ecurrence pr\'ec\'edente. $\hfill \blacksquare$

\section{Calculs d'exemples.}

\subsection{Le cas de \ $x^5 + y^5 + x^2.y^2.$}

\subsubsection{Annuler le g\'en\'erateur.}

Dans ce paragraphe nous commencerons par  expliciter un polyn\^ome en (a,b) qui annule l'\'el\'ement \ $1$ \ du (a,b)-module de Brieskorn de \ $f(x,y) : = x^5 + y^5 + x^2.y^2.$ Ceci correspond aux int\'egrales du type
$$ \int_{\gamma_s} \frac{dx\wedge dy}{df}  $$
o\`u \ $(\gamma_s)_{s \in D}$  \ est une famille horizontale de \ $1-$cycles compacts dans les fibres de \ $f$

\bigskip

\begin{lemma}\label{1}
On a
$$ x^{10}.y^{10} = (5a - 18b)(5a - 14b)(5a - 10b)(5a - 6b)(5a - 2b)(1) .$$
\end{lemma}

\parag{Preuve} Commen{\c c}ons par montrer que l'on a
\begin{equation*}
 x^{p+5}.y^p = x^p.y^{p+5} \tag{a}
 \end{equation*}
pour tout entier \ $p \geq 0$. On a
\begin{align*}
& b(x^p.y^p) = (5x^4 + 2x.y^2)\frac{1}{p+1}x^{p+1}.y^p = \frac{1}{p+1}\big[5x^{p+5}.y^p + 2x^{p+2}.y^{p+2}\big] \tag{b} \\
& \qquad       =  (5y^4 + 2x^2.y)\frac{1}{p+1}x^{p}.y^{p+1} = \frac{1}{p+1}\big[5x^{p}.y^{p+5} + 2x^{p+2}.y^{p+2}\big]
\end{align*}
ce qui donne l'\'egalit\'e annonc\'ee. On en d\'eduit que
\begin{equation*}
a(x^p.y^p) = 2x^{p+5}.y^p + x^{p+2}.y^{p+2} \tag{c} 
\end{equation*}
et en combinant avec \ $(b)$ \ cela donne
\begin{equation*}
(5a - 2(p+1).b)(x^p.y^p) = x^{p+2}.y^{p+2} \tag{d} 
\end{equation*}
ce qui permet de conclure. $\hfill \blacksquare$

\begin{lemma}
On a
\begin{equation*}
(a - 4b)(x^5.y^{10}) = \frac{-1}{2}.x^{10}.y^{10}  \tag{e}
\end{equation*}
\end{lemma}

\parag{Preuve}
On a
\begin{align*}
& b(x^5.y^{10}) = (5x^4 + 2x.y^2)\frac{1}{6}x^6.y^{10} = \frac{1}{6}\big[5x^{10}.y^{10} + 2x^7.y^{12}\big] \tag{f} \\
& b(x^5.y^{10}) = (5y^4 + 2x^2.y)\frac{1}{11}x^5.y^{11} = \frac{1}{11}\big[5x^5.y^{15} + 2x^7.y^{12}\big] \tag{g} \\
& a(x^5.y^{10}) =  x^{10}.y^{10} + x^5.y^{15} + x^7.y^{12} \tag{h}
\end{align*}
ce qui donne successivement, en combinant \ $(h)$ \ et \ $(g)$ \ puis  avec \ $(f)$
\begin{align*}
& (5a - 11b)(x^5.y^{10})  = 5x^{10}.y^{10} + 3x^7.y^{12} \\
& (10a - 40b)(x^5.y^{10}) = -5x^{10}.y^{10} 
\end{align*}
d'o\`u l'\'egalit\'e annonc\'ee. $\hfill \blacksquare$

\bigskip

\begin{lemma}
On a aussi
\begin{align*}
& (2a - b)(1) = -x^5 \tag{i}\\
& (2a - 3b)(x^5) = -x^5.y^5 \tag{i'} \\
& (a - 3b)(x^5.y^5) = \frac{-1}{2}x^5.y^{10}  \tag{i''} 
\end{align*}
\end{lemma}

\parag{Preuve} Comme on a \ $x^5 = y^5$, on aura \ $a(1) = 2x^5 + x^2.y^2 $ \ et \ $b(1) = 5x^5 + 2x^2.y^2$. D'o\`u la premi\`ere \'egalit\'e. De m\^eme, \ $a(x^5) = x^{10} + x^5.y^5 + x^7.y^2$ \ et \ $$b(x^5) = (5x^4 + 2x.y^2)\frac{x^6}{6} = \frac{1}{6}[5x^{10} + 2x^7.y^2] = (5y^4 + 2x^2.y)x^5.y = 5x^5.y^5 + 2x^7.y^2 $$
ce qui donne \ $ x^{10} = 6x^5.y^5 + 2x^7.y^2 $ \ d'o\`u \ $a(x^5) =  7x^5.y^5 + 3x^7.y^2$ \ et permet de conclure pour \ $\rm{(i')}$.\\
Utilisons \ $(a)$ \ pour \ $p = 5$ :
\begin{align*}
& b(x^5.y^5) = (5x^4 + 2xy^2)\frac{1}{6}x^6.y^5 =  \frac{1}{6}\big[5x^{10}.y^{5} + 2x^7.y^{7}\big] \quad {\rm et}  \tag{k}\\
& a(x^5.y^5) = 2x^{10}.y^5 + x^7.y^7 \quad {\rm donnent} \\
& (a - 3b)(x^5.y^5) = \frac{-1}{2}x^5.y^{10} 
\end{align*}
et donc \ $\rm{(i'')}$. $\hfill \blacksquare$

\parag{Conclusion}
On a donc
$$ 4(a - 4b)(a - 3b)(2a - 3b)(2a - b)(1) =  (5a - 18b)(5a - 14b)(5a - 10b)(5a - 6b)(5a - 2b)(1) $$
et donc l'\'el\'ement
$$  (5a - 18b)(5a - 14b)(5a - 10b)(5a - 6b)(5a - 2b) -  (2a - 8b)(2a - 6b)(2a - 3b)(2a - b) $$
est dans l'annulateur de \ $1$ \ dans le module de Brieskorn de \ $f$.

\bigskip

Traitons maintenant le cas du (a,b)-module engendr\'e par le mon\^ome \ $x$.

\begin{lemma} On a 
$$ x^{11}.y^{10} =  (5a -19b)(5a - 15b)(5a - 11b)(5a - 7b)(5a - 3b)[x] .$$
\end{lemma}

\parag{Preuve} Posons \ $\alpha_p : = x^{2p+1}.y^{2p} $ \ pour chaque \ $p \geq 0$. On a
\begin{align*}
& b(\alpha_p) = (5x^4 + 2x.y^2)\frac{x^{2p+2}.y^{2p}}{2p+2} = \frac{1}{2p+2}\big[5x^{2p+6}.y^{2p}+ 2x^{2p+3}.y^{2p+2}\big] \\
&  b(\alpha_p) = (5y^4 + 2x^2.y)\frac{x^{2p+1}.y^{2p+1}}{2p+1} = \frac{1}{2p+1}\big[5x^{2p+1}.y^{2p+5} + 2x^{2p+3}.y^{2p+2}\big] \\
& a(\alpha_p) = x^{2p+6}.y^{2p} + x^{2p+1}.y^{2p+5} + x^{2p+3}.y^{2p+2}
\end{align*}
On en d\'eduit que l'on a pour chaque \ $p \geq 0$
\begin{equation*}
(5a - (4p+3)b)[\alpha_p] = \alpha_{p+1} 
\end{equation*}
ce qui permet de conclure puisque \ $\alpha_0 = x$ \ et \ $\alpha_5 = x^{11}.y^{10}$. $\hfill \blacksquare$

 \begin{lemma}
  On a 
   $$ (a - \frac{42}{10}b)(a - \frac{27}{10}b)(a - \frac{23}{10}b)(a - \frac{8}{10}b)(x) = \frac{1}{16}.x^{11}.y^{10}   .$$
  \end{lemma}
  
  \parag{Preuve} Calculons successivement les mon\^omes \ $x^6, x^{11}, x^{11}.y^5$ \ et \ $x^{11}.y^{10}$.
  \begin{align*}
  &  b(x) = (5x^4 + 2x.y^2)\frac{x^2}{2} = \frac{1}{2}\big[5x^6 + 2x^3.y^2\big] 
  \end{align*}
  ce qui donne, puisque \ $x^3.y^2 = \alpha_1 = (5a - 3b)(x)$
  \begin{align*}
  & 5x^6 = 2b(x) - 2(5a - 3b)(x) = (5b - 10a)(x) \quad {\rm soit} \\
  & -5x^6 =  (10a - 8b)(x) \tag{t}
  \end{align*}
  Puis
  \begin{align*}
  & b(x^6) = (5x^4 + 2x.y^2)\frac{x^7}{7} = \frac{1}{7}\big[5x^{11} + 2x^8.y^2\big] \\
  & b(x^6) = (5y^4 + 2x^2.y)x^6.y = 5x^6.y^5 + 2x^8.y^2 \\
  & a(x^6) = x^{11} + x^6.y^5 + x^8.y^2
  \end{align*}
  ce qui donne
  \begin{equation*}
  -5x^{11} = (10a - 23b)(x^6)   \tag{u}
  \end{equation*}
  Ensuite
  \begin{align*}
  & b(x^{11}) = (5x^4 + 2x.y^2)\frac{x^{12}}{12} = \frac{1}{12}\big[5x^{16} + 2x^{13}.y^2\big] \\
   & b(x^{11}) = (5y^4 + 2x^2.y)x^{11}.y = 5x^{11}.y^5 + 2x^{13}.y^2 \\
   & a(x^{11}) = x^{16} + x^{11}.y^5 + x^{13}.y^2 
   \end{align*}
   ce qui donne
   \begin{equation*}
   -5x^{11}.y^5 = (10a - 27b)(x^{11}) \tag{v}
   \end{equation*}
   Et enfin
   \begin{align*}
   & b(x^{11}.y^5) = (5x^4 + 2x.y^2)\frac{x^{12}.y^5}{12} = \frac{1}{12}\big[5x^{16}.y^5 + 2x^{13}.y^7\big] \\
   & b(x^{11}.y^5) = (5y^4 + 2x^2.y)\frac{x^{11}.y^6}{6} =  \frac{1}{6}\big[5x^{11}.y^{10} + 2x^{13}.y^7\big] \\
    & a(x^{11}.y^5) = x^{16}.y5 + x^{11}.y^{10} + x^{13}.y^7
    \end{align*}
    ce qui donne
    \begin{equation*}
    - 5x^{11}.y^{10} = (10a - 42b)(x^{11}.y^5) \tag{w}
    \end{equation*}
    On en conclut donc que l'on a bien
    $$ (10a - 42b)(10a - 27b)(10a -23b)(10a - 8b)(x) = 625.x^{11}.y^{10} .\qquad \qquad \blacksquare $$
    
    \bigskip
    
    Traitons enfin le cas du mon\^ome \ $xy$.
    
    \begin{lemma} Posons \ $\alpha_p : = x^{2p+1}.y^{2p+1} $ \ pour \ $p \geq 0$ . On a 
    $$ (5a - 4(p+1)b)(\alpha_p) = \alpha_{p+1}  $$
    ce qui donne
    $$ x^{11}.y^{11} = (5a - 20b)(5a - 16b)(5a -12b)(5a -8b)(5a - 4b)[xy] .$$
    \end{lemma}
    
   \noindent  Le calcul,  analogue aux pr\'ec\'edents, est laiss\'e au lecteur.
    
    \begin{lemma}
    On a
    $$ \frac{1}{16}x^{11}.y^{11} = (2a - 9b)(2a - 6b)(2a - 5b)(2a - b)[xy] .$$
    \end{lemma}
    
    \parag{Preuve} Le lecteur, en s'inspirant des calculs pr\'ec\'edents,  prouvera sans difficult\'es les relations suivantes :
    \begin{align*}
    & (2a - 2b)[xy] = - x^6.y \\
    & (2a - 5b)(x^6.y) = - x^{11}.y \\
    & (2a - 6b)(x^{11}.y) = - x^{11}.y^6 \\
    & (2a - 9b)(x^{11}.y^6) = - x^{11}.y^{11} 
    \end{align*}
    ce qui prouve notre assertion. $\hfill \blacksquare$

\subsubsection{Application}

 Le but de ce paragraphe est de montrer que le (a,b)-module (monog\`ene) engendr\'e par \ $1$ \ dans le module de Brieskorn de \ $f(x,y) = x^5 + y^5 + x^2.y^2$ \ est de rang \ $4$ \ et que  son polyn\^ome de Bernstein est 
 $$(x+1)^2.(x+\frac{1}{2})^2.$$
 
\noindent  Pour le mon\^ome \ $x$ \ on trouve  que le polyn\^ome de Bernstein vaut
 $$ (x + \frac{6}{5})(x + \frac{4}{5})(x + \frac{13}{10})(x + \frac{7}{10}) .$$
 
 \noindent Pour le mon\^ome \ $xy$ \ on trouve que le polyn\^ome de Bernstein vaut
$$ (x + \frac{3}{2})^2(x + 1)(x + \frac{1}{2}) . $$
 
  \bigskip
 
\noindent Pour cela, compte tenu des paragraphes pr\'ec\'edents, il nous suffit de prouver que les \'el\'ements \ $1, a(1), a^2(1),a^3(1)$ \ sont lin\'eairement ind\'ependants dans le \ $\mathbb{C}-$espace vectoriel \ $\tilde{\mathcal{A}}.1\big/b.\tilde{\mathcal{A}}.1$.\\
La v\'erification analogue pour les mon\^omes \ $x$ \ et \ $xy$ \ est laiss\'ee au lecteur.\\

\noindent  Considerons donc une relation du type
 $$ \lambda_0.1 + \lambda_1.a(1) + \lambda_2.a^2(1) + \lambda_3.a^3(1) \in b.\tilde{\mathcal{A}}.1.$$
 Comme \ $1$ \ et \ $x^2.y^2$ \ sont lin\'eairement ind\'ependants modulo l'id\'eal jacobien de \ $f$, on en d\'eduit que n\'ecessairement on a \ $\lambda_0 = \lambda_1 = 0$.\\
  Montrons que \ $x^4.y^4$ \ et \ $x^6.y^6$ \ ne sont pas dans \ $ b.\tilde{\mathcal{A}}.1$. Cela r\'esulte de la liste suivante\footnote{on utilise ici les \'egalit\'es de sym\'etrie en (x,y) prouv\'ees au lemme \ref{1}.} qui engendre les polyn\^omes  de degr\'es  \ $\leq 15$ \ en \ $(x,y)$ \ qui sont dans \ $ b.\tilde{\mathcal{A}}.1$:
  \begin{align*}
  & b(1) = 5x^5 + 2x^2y^2 \\
  & b(x^2.y^2) = \frac{5}{2}x^7.y^2 + x^4.y^4 \\
  & b(x^5) = \frac{1}{6}[5x^{10} + 2x^7.y^2] \\
  & b(x^4.y^4) = x^9.y^4 + \frac{2}{5}x^6.y^6 \\
  & b(x^7.y^2) = \frac{1}{8}[x^{12}.y^2 + 2x^9.y^4] \\
  & b(x^{10}) = \frac{1}{11}[5x^{15} + 2x^{12}.y^2]
  \end{align*}

 \bigskip
 
  Comme \ $25a^2(1) = x^4.y^4 \  modulo \ \ b.\tilde{\mathcal{A}}.1$, et comme \ $x^6.y^6 = 5a(x^4.y^4) \  modulo \ \ b.\tilde{\mathcal{A}}.1$, on obtient la relation
   $$5\lambda_2.x^4.y^4 + \lambda_3.a(x^4.y^4) \in b.\tilde{\mathcal{A}}.1.$$
   Comme on a vu que \ $\lambda_2$ \ et \ $\lambda_3$ \ sont non nuls, ou alors tous les deux nuls, leur non nullit\'e donnerait une relation du type \ $a(x^4.y^4) - \rho.x^4.y^4 \in b.\tilde{\mathcal{A}}.1$ \ avec \ $\rho \not= 0 $. Ce qui donnerait \ $a^n(x^4.y^4) - \rho^n.x^4.y^4 \in b.\tilde{\mathcal{A}}.1$ \ pour tout \ $n \in \mathbb{N}$. Comme \ $a^n.\tilde{\mathcal{A}}.1\subset b.\tilde{\mathcal{A}}.1$ \ pour \ $n$ \ assez grand (par r\'egularit\'e), on obtiendrait \ $x^4.y^4 \in b.\tilde{\mathcal{A}}.1$. ce qui est absurde. Donc \ $\lambda_2 = \lambda_3 = 0$ \ et le rang de \ $\tilde{\mathcal{A}}.1$ \ est \ $\geq 4$.\\
   Le paragraphe pr\'ec\'edent permet alors de conclure gr\^ace au th\'eor\`eme \ref{Utile}. $\hfill \blacksquare$

\subsection{Le cas de \ $f : = x^n + y^n + z^n + x.y.z $.}

On suppose \ $n \geq 4$, le cas \ $n = 3$ \ qui est homog\`ene est imm\'ediat.

\begin{lemma}
Pour tout \ $p \geq 0 $ \ on a
$$ (a - \frac{3(p+1)}{n}b)(x^p.y^p.z^p) = \frac{n-3}{n} x^{p+1}.y^{p+1}.z^{p+1} .$$
\end{lemma}

\parag{Preuve} On a
\begin{align*}
& b(x^p.y^p.z^p) = (nx^{n-1} + y.z)\frac{x^{p+1}}{p+1}.y^p.z^p = \frac{1}{p+1}[nx^{n+p}.y^p.z^p + x^{p+1}.y^{p+1}.z^{p+1}]
\end{align*}
ce qui donne \ $ x^{n+p}.y^p.z^p = y^{n+p}.z^p.x^p = etc ... $ \ et donc
\begin{align*}
& a(x^p.y^p.z^p) = 3x^{n+p}.y^p.z^p + x^{p+1}.y^{p+1}.z^{p+1} 
\end{align*}
d'o\`u la relation annonc\'ee. $\hfill \blacksquare$

On notera que l'on a montr\'e au passage que l'on a \ $x^n = y^n = z^n$ \ dans le (a,b)-module de Brieskorn de \ $f$.

\begin{lemma}
On a 
$$ (a - 2b)(x^n) = (3-n)x^n.y^n .$$
\end{lemma}

\parag{Preuve} Puisque l'on a \ $x^n = y^n$, cela donne :
\begin{align*}
& b(x^n) = (nx^{n-1} + y.z)y^n.x = nx^n.y^n + x.y^{n+1}.z \\
&  b(x^n) = (nx^{n-1} + y.z)\frac{x^{n+1}}{n+1} = \frac{1}{n+1}[nx^{2n} + x^{n+1}.y.z ]
\end{align*}
ce qui permet d'obtenir \ $ nx^{2n} = (n+1)b(x^n) - x^{n+1}.y.z $ \ et donc, en utilisant le fait que  \ $x^n.y^n = x^n.z^n = y^n.z^n  : = v$ \ et en posant \ $x^{n+1}.y.z = x.y^{n+1}.z = x.y.z^{n+1} : = u$
\begin{align*}
& na(x^n) =   (n+1)b(x^n) - u + 2v + u \quad {\rm et \ donc} \\
&  a(x^n) = b(x^n) + 3v + u = 2b(x^n) + (3-n)v
\end{align*}
puisque \ $b(x^n) = nv + u $. $\hfill \blacksquare$

Calculons enfin \ $b(x^n.y^n)$ :
\begin{align*}
& b(x^n.y^n) = (nz^{n-1} + x.y)x^n.y^n.z = nx^n.y^n.z^n + x^{n+1}.y^{n+1},z \\
& \quad          = (nx^{n-1} + yz)\frac{x^{n+1}}{n+1}y^n = \frac{1}{n+1}[nx^{2n}.y^n + x^{n+1}.y^{n+1}.z ]
\end{align*}
On en d\'eduit les \'egalit\'es
\begin{align*}
& x^{n+1}.y^{n+1}.z = z^{n+1}.y^{n+1}.z = x^{n+1}.z^{n+1}.y : = s \\
& x^{2n}.y^n = y^{2n}.x^n = \cdots : = t 
\end{align*}
En posant \ $x^n.y^n.y^n : = r$ \ cela donne
\begin{align*}
& b(x^n.y^n) = nr + s = \frac{1}{n+1}[nt + s] \\
& a(x^n.y^n) = t + t + r + s 
\end{align*}
On en d\'eduit la relation
\begin{align*}
(a - 3b)(x^n,y^n) = (3-n)x^n.y^n.z^n 
\end{align*}
Comme le premier lemme donne
$$\big(\prod_{p=0}^{n-1} (a - \frac{3(p+1)}{n}b)\big)(1) = (\frac{n-3}{n})^n.x^n.y^n.z^n $$
on obtient que
$$ \frac{n^n}{(n-3)^{n-3}}.\big(\prod_{p=0}^{n-1} (a - \frac{3(p+1)}{n}b)\big) + (a-3b)(a-2b)(a-b) $$
annule \ $1$ \ dans le module de Brieskorn de \ $f$.\\
On en conclut que le polyn\^ome de Bernstein de \ $\A.1$ \ est \'egal \`a \ $(x + 1)^3$.

\subsection{Commentaire final.}

Si on consid\`ere un (a,b)-module de la forme \ $\A\big/\A.(a^4 + a^5)$ \ (analogue au cas du paragraphe 5.1), on constate que l'on a une suite exacte
$$ 0 \to K \to \A\big/\A.(a^4 + a^5) \to \A\big/\A.a^4 \to 0 $$
et le quotient \ $\A/\A.a^4 $ \ est monog\`ene r\'egulier de rang \ $4$. Mais le noyau est de rang \ $1$ \ et non r\'egulier : en effet le noyau est \ $K \simeq \A\big/\A.(a + 1) $ \ n'est m\^eme pas local (donc encore moins r\'egulier).\\
Quand on tensorise cette suite exacte par le compl\'et\'e \ $a-$adique \ $\hat{A}$ \  de \ $\A$, on a \ $K \otimes \hat{A} = 0$ \ et il ne reste que \ $\A\big/\A.a^4 \otimes \hat{A}  \simeq \hat{A}\big/\hat{A}.a^4 $.\\
On notera qu'un (a,b)-module local (donc \`a fortiori s'il est r\'egulier)  \ $E$ \ v\'erifie toujours
$$ E \simeq E \otimes_{\A} \hat{A} $$
puisqu'il est complet pour la filtration \ $a-$adique\footnote{O utilise ici le fait que \ $\hat{A}$ \ est un \ $\A-$module \`a gauche et \`a droite.}. \\
Ceci explique pourquoi dans les exemples pr\'ec\'edents le rang ne correspond pas au degr\'e en \ $a$ \ des polyn\^omes trouv\'es.\\
Mais par exemple \ $ a = -1$ \ qui correspond au noyau \ $K$ \ exprime le fait que l'\'equation diff\'erentielle satisfaite par l'int\'egrale consid\'er\'ee admet un (autre) point singulier (ici  le point \ $-1$). Cela sugg\`ere de consid\'erer les diff\'erentes localisations int\'eressantes d'un (a,b)-module (alg\'ebrique) et d'avoir une th\'eorie "globale" sur \ $\mathbb{C}$ \ des \'equations diff\'erentielles alg\'ebriques. La localisation au point \ $z$ \ signifie que l'on consid\`ere sur \ $E$ \ les op\'erateurs \ $a_z : = a + z $ \ et \ $b$ \ et que l'on compl\`ete pour la topologie \ $a_z-$adique. Cela ressemble fort \`a la th\'erie des schemas ! Il est clair que les int\'egrales de p\'eriodes de fonctions alg\'ebriques  (ou de vari\'et\'es alg\'ebriques) doivent rentrer dans un cadre de ce genre.

\bigskip

\newpage

\section{R\'ef\'erences.}

\begin{itemize}

 \item{[Br.70]} Brieskorn, E. {\it Die Monodromie der Isolierten Singularit{\"a}ten von Hyperfl{\"a}chen}, Manuscripta Math. 2 (1970), p. 103-161.

\item{[B. 93]} Barlet, D. {\it Th\'eorie des (a,b)-modules I}, in Complex Analysis and Geometry, Plenum Press, (1993), p. 1-43.

\item{[B. 95]} Barlet, D. {\it Th\'eorie des (a,b)-modules II. Extensions}, in Complex Analysis and Geometry, Pitman Research Notes in Mathematics Series 366 Logman (1997), p. 19-59.

\item{[B. 05]} Barlet, D. {\it Module de Brieskorn et forme hermitiennes pour une singularit\'e isol\'ee d'hypersuface}, revue de l'Inst. E. Cartan (Nancy) 18 (2005), p. 19-46.

\item{[B. II]} Barlet, D. {\it Sur certaines singularit\'es d'hypersurfaces II}, J. Alg. Geom. 17 (2008), p. 199-254.

\item{[B. 07]} Barlet, D. {\it Sur les fonctions a singularit\'e de dimension 1(version r\'evis\'ee)}, preprint Institut E. Cartan (Nancy) (2008) $n^042$, p. 1-26.

\item{[B. 08]} Barlet, D. {\it Two finiteness theorem for regular (a,b)-modules}, preprint Institut E. Cartan (Nancy) (2008) $n^0 5$, p. 1-38.

\item{[B.-S. 04]} Barlet, D. et Saito, M. {\it Brieskorn modules and Gauss-Manin systems for non isolated hypersurface singularities,} J. Lond. Math. Soc. (2) 76 (2007) $n^01$ \ p. 211-224.

\item{[M. 75]} Malgrange, B. {\it Le polyn\^ome de Bernstein d'une singularit\'e isol\'ee}, in Lect. Notes in Math. 459, Springer (1975), p.98-119.

\item{[S. 89]} Saito, M. {\it On the structure of Brieskorn lattices}, Ann. Inst. Fourier 39 (1989), p.27-72.

\end{itemize}

\end{document}